\documentclass[11pt]{article}
\usepackage{amsfonts,amscd,amssymb}

\input{diagrams}
\newtheorem{definition}{Definition}[section]
\newtheorem{lemma}[definition]{Lemma}
\newtheorem{proposition}[definition]{Proposition}

\newtheorem{remark}[definition]{Remark}

\newtheorem{theorem}[definition]{Theorem}
\newtheorem{example}[definition]{Example}

\newcommand{\Ent}{\mathbb{Z}}

\newcommand{\gsm}{\mbox{$\blacktriangleright \hspace{-0.7mm}<$}}

\def\rawo\lonra{\longrightarrow}

\def\ot{\otimes}

\newcommand{\eqref}[1]{(\ref{eq:#1})}

\newenvironment{proof}{{\it Proof.}}{\hfill $ \square $ \vskip 4mm}
\begin{document}
\title{Extending lazy 2-cocycles on Hopf algebras and\\
lifting projective representations afforded by them }
\author{Juan Cuadra\thanks{The first author was partially supported by the
DGES-project BMF2002-02717.} \\
Dept. \'Algebra y An\'alisis Matem\'atico\\
Universidad de Almer\'{\i}a,
E-04120 Almer\'{\i}a, Spain \\
e-mail: jcdiaz@ual.es \and Florin Panaite\thanks {Research
carried out while the second author was visiting the University of
Almer\'{\i}a supported by a NATO fellowship offered by the Spanish
Ministry of Science and Technology. This author was also
partially supported by the programme CERES of the Romanian
Ministry of Education and
Research, contract no. 4-147/2004.   }\\
Institute of Mathematics of the
Romanian Academy\\
PO-Box 1-764, RO-014700 Bucharest, Romania\\
 e-mail: Florin.Panaite@imar.ro}
\date{}
\maketitle

\begin{abstract}
We study some problems related to lazy 2-cocycles, such as: extension of
(lazy) 2-cocycles to a Drinfeld double and to a Radford biproduct,
Yetter-Drinfeld data obtained from lazy 2-cocycles, lifting of
projective representations afforded by lazy 2-cocycles.
\end{abstract}

\section*{Introduction}

A left 2-cocycle $\sigma :H\ot H\rightarrow k$ on a Hopf algebra
$H$ is called {\it lazy} if it satisfies the condition
\begin{eqnarray*}
&&\sigma (h_1, h'_1)h_2h'_2=\sigma (h_2, h'_2)h_1h'_1, \;\;\;\;\;
\forall \;h, h'\in H.
\end{eqnarray*}
This kind of cocycles were used in \cite{gio1} as a tool to
compare the Brauer groups of Sweedler's Hopf algebra with respect
to the different quasitriangular structures. See also \cite{cc1}
and \cite{cc} for an application of this technique to other sort
of Hopf algebras. Lazy cocycles and lazy cohomology were also
used in \cite{sch} to give a generalized version of Kac's exact
sequence. A general theory of lazy cocycles and lazy cohomology
started to be developed recently in \cite{bc}. The most
remarkable fact is that the set $Z^2_L(H)$ of normalized and
convolution invertible lazy 2-cocycles on $H$ form a group, and
that one can also define lazy 2-coboundaries $B^2_L(H)$ and the
second lazy cohomology group $H^2_L(H)=Z^2_L(H)/B^2_L(H)$,
generalizing the second Sweedler cohomology group of a
cocommutative Hopf algebra (note that for cocommutative Hopf
algebras {\it any} 2-cocycle is lazy). The group $H^2_L(H)$ can
be regarded as a subgroup of $Bigal(H)$, the group of Bigalois
objects of $H$, and the examples in \cite{bc} show that it is
much easier to compute $H^2_L(H)$ than $Bigal (H)$. \par
\smallskip

In general, the results in \cite{bc} suggest that, for an
arbitrary Hopf algebra, lazy cocycles are much closer to the
cocommutative case than general left cocycles. Hence, a sort of
general principle is suggested: results that hold for an {\it
arbitrary} 2-cocycle on a {\it cocommutative} Hopf algebra are
likely to hold also for a {\it lazy} 2-cocycle on an {\it
arbitrary} Hopf algebra. A good example of this principle is the
extension of Schur-Yamazaki formula in \cite{bc} that allows to
describe the second lazy cohomology group of a tensor product of
Hopf algebras. Throughout this paper we will verify this
principle several times. \par \smallskip

This paper is a contribution to the study of lazy cocycles and
lazy cohomology, in three different directions: the problem of
extending (lazy) 2-cocycles to a Drinfeld double and to a Radford
biproduct; Yetter-Drinfeld data obtained from lazy 2-cocycles;
lifting of projective representations afforded by lazy
2-cocycles. As we will see below, each of these directions has
its own (natural) motivations and possible applications.\par
\smallskip

We describe now in some detail the contents of the paper. After
presenting in Section \ref{sec1} some preliminaries, in Section \ref{sec2} we
provide some new properties of lazy 2-cocycles that are needed in
the next sections, but which could also be of independent
interest. Among these properties is the following formula:
\begin{eqnarray*}
&&\sigma (h_1, S(h_2))=\sigma (S(h_1), h_2), \;\;\;\;\;\forall \;h\in H,
\end{eqnarray*}
for a lazy 2-cocycle $\sigma$ on a Hopf algebra $H$; this formula
is important and well-known for group algebras, and we show that
in general it is false if $\sigma $ is not lazy. \par \smallskip

In Section \ref{sec3} we prove that any lazy 2-cocycle $\sigma$ on a
finite dimensional Hopf algebra $H$ can be extended to a lazy
2-cocycle $\overline{\sigma}$ on the Drinfeld double $D(H)$ (this
property can be obtained also from results in \cite{bc}, where
moreover a complete description of $H^2_L(D(H))$ is given). We
point out that this extension is canonical in a certain sense
(expressed in terms of the so-called {\it diagonal crossed
product}, a construction introduced in \cite{hn1}; actually, the
relation between lazy 2-cocycles and the diagonal crossed product
was our starting point for this article). Section \ref{sec4} is devoted to
the extension of cocycles on a Radford biproduct. We consider a
Radford biproduct $B\times H$, with $H$ a Hopf algebra and $B$ a
Hopf algebra in the Yetter-Drinfeld category $_H^H{\cal YD}$.
Cocycles and the second lazy cohomology group $H^2_L(B)$ may be
defined in the category $_H^H{\cal YD}$. We prove that, if
$\sigma $ is a left 2-cocycle on $B$ in $_H^H{\cal YD}$, it can
be extended canonically to a left 2-cocycle $\overline{\sigma}$
on $B\times H$, $\sigma $ lazy in $_H^H{\cal YD}$ implies
$\overline{\sigma }$ lazy and the map $\sigma \mapsto
\overline{\sigma }$ induces a group morphism $H^2_L(B)\rightarrow
H^2_L(B\times H)$. \par \smallskip

In Section \ref{sec5} we study Yetter-Drinfeld data obtained from lazy
2-cocycles. Namely, if $\sigma :H\ot H\rightarrow k$ is a
normalized and convolution invertible lazy 2-cocycle, we have the
$H$-bicomodule algebra $H(\sigma )=$$\;_{\sigma }H=H_{\sigma }$,
hence the Yetter-Drinfeld category $_{H(\sigma )}{\cal YD}(H)^H$.
We prove that, if $M$ is a finite dimensional object in this
category, then $End(M)$ and $End(M)^{op}$ are algebras in
$_H{\cal YD}^H$ (we hoped that they were even Azumaya in $_H{\cal
YD}^H$, but in general they are not). More can be said if $H$ is
finite dimensional (for this we use again the diagonal crossed
product and results from Section \ref{sec3}). This section is partially
motivated by the belief (expressed also in \cite{bc}) that it
should exist a relation between $H^2_L(H)$ and the Brauer group
$BQ(k, H)$ of $H$ (hopefully, an embedding of $H^2_L(H)$ into
$BQ(k, H)$), at least for some classes of Hopf algebras (for
instance, the cotriangular ones). Finally, in Section \ref{sec6} we prove
that any Hopf algebra $H$ admits a central extension $B$ with the
property that any projective representation of $H$ afforded by a
{\it lazy} 2-cocycle can be lifted to an ordinary representation
of $B$. The case when $H$ is cocommutative was done by I. Boca
(generalizing in turn the classical case of groups, due to
Schur); our proof follows closely the one of Boca. This section
could be regarded as a good illustration of the general principle
we mentioned before (replacement of cocommutativity by laziness).

\setcounter{section}{0}
\section{Preliminaries}\label{sec1}

\setcounter{equation}{0}

In this section we recall some definitions and results and we fix
some notation that will be used throughout the paper. \par
\smallskip

We will work over a ground field $k$. All algebras, linear
spaces, etc, will be over $k$; unadorned $\ot$ means $\ot _k$.
For a Hopf algebra $H$ with comultiplication $\Delta$ we use the
version of Sweedler's sigma notation: $\Delta (h)=h_1\ot h_2$.
Unless otherwise stated, $H$ will denote a Hopf algebra with
bijective antipode $S$. For a linear map $\sigma :H\ot
H\rightarrow k$ we will use either the notation $\sigma (h,h')$
or $\sigma(h\ot h')$. \par \smallskip

A linear map $\sigma :H\ot H\rightarrow k$ is called a {\it left
2-cocycle} if it satisfies the condition
$$\sigma (a_1, b_1)\sigma (a_2b_2, c)=\sigma (b_1, c_1)\sigma (a,
b_2c_2),$$ for all $a, b, c\in H$, and it is called a {\it right
2-cocycle} if it satisfies the condition $$\sigma (a_1b_1,
c)\sigma (a_2, b_2)=\sigma (a, b_1c_1)\sigma (b_2, c_2).$$ Given
a linear map $\sigma :H\ot H\rightarrow k$, define a product
$\cdot_{\sigma }$ on $H$ by
\begin{eqnarray*}
&&h\cdot _{\sigma }h'=
\sigma (h_1, h'_1)h_2h'_2, \;\;\;\;\;\forall \;h, h'\in H.
\end{eqnarray*}
Then $\cdot_{\sigma }$ is associative if and only if $\sigma $ is
a left 2-cocycle. If we define $\cdot _{\sigma }$ by
\begin{eqnarray*}
&&h\cdot _{\sigma }h'=h_1h'_1\sigma (h_2, h'_2), \;\;\;\;\;\forall \;
h, h'\in H,
\end{eqnarray*}
then $\cdot _{\sigma }$ is associative if and only if $\sigma $
is a right 2-cocycle. In any of the two cases, $\sigma$ is
normalized (i.e. $\sigma (1, h)= \sigma (h, 1)=\varepsilon (h)$
for all $h\in H$) if and only if $1_H$ is the unit for $\cdot
_{\sigma }$.  If $\sigma $ is a normalized left (respectively
right) 2-cocycle, we denote the algebra $(H, \cdot _{\sigma })$
by $_{\sigma }H$ (respectively $H_{\sigma }$). It is well-known
that $_{\sigma }H$ (respectively $H_{\sigma }$) is a right
(respectively left) $H$-comodule algebra via the comultiplication
$\Delta $ of $H$. If $\sigma :H\ot H\rightarrow k$ is normalized
and convolution invertible, then $\sigma $ is a left 2-cocycle if
and only if $\sigma ^{-1}$ is a right 2-cocycle. \par \smallskip

If $\gamma :H\rightarrow k$ is linear, normalized (i.e. $\gamma
(1)=1$) and convolution invertible, define $D^1(\gamma ):H\ot
H\rightarrow k$ by
\begin{eqnarray*}
&&D^1(\gamma )(h, h')=\gamma (h_1)\gamma (h'_1)\gamma ^{-1}(h_2h'_2),
\;\;\;\;\;\forall \;h, h'\in H.
\end{eqnarray*}
Then $D^1(\gamma )$ is a normalized and convolution invertible
left 2-cocycle. If $\sigma, \sigma':H\ot H\rightarrow k$ are
normalized and convolution invertible left 2-cocycles, they are
called {\it cohomologous} if there exists $\gamma :H\rightarrow k$
normalized and convolution invertible such that
\begin{eqnarray*}
&&\sigma '(h, h')=\gamma (h_1)\gamma (h'_1)\sigma (h_2, h'_2)
\gamma ^{-1}(h_3h'_3), \;\;\;\;\;\forall \;h, h'\in H.
\end{eqnarray*}

We recall now from \cite{bc} some facts about lazy cocycles and
lazy cohomology. The set $Reg^1 (H)$ (respectively $Reg^2 (H)$)
consisting of normalized and convolution invertible linear maps
$\gamma :H\rightarrow k$ (respectively $\sigma :H\ot H\rightarrow
k$), is a group under the convolution product. An element $\gamma
\in Reg^1 (H)$ is called {\it lazy} if
\begin{eqnarray*}
&&\gamma (h_1)h_2=h_1\gamma (h_2), \;\;\;\;\;\forall \;h\in H.
\end{eqnarray*}
The set of lazy elements of $Reg^1 (H)$, denoted by $Reg^1_L (H)$,
is a central subgroup of $Reg^1 (H)$. An element $\sigma \in
Reg^2 (H)$ is called {\it lazy} if
\begin{eqnarray*}
&&\sigma (h_1, h'_1)h_2h'_2=h_1h'_1\sigma (h_2, h'_2),\;\;\;\;\;
\forall \;h, h'\in H.
\end{eqnarray*}
The set of lazy elements of $Reg^2 (H)$, denoted by $Reg^2_L (H)$,
is a subgroup of $Reg^2 (H)$. We denote by $Z^2 (H)$ the set of
left 2-cocycles on $H$ and by $Z^2_L (H)$ the set $Z^2 (H)\cap
Reg^2_L (H)$ of normalized and convolution invertible lazy
2-cocycles. If $\sigma \in Z^2_L(H)$, then the algebras $_{\sigma
}H$ and $H_{\sigma }$ coincide and will be denoted by $H(\sigma
)$; moreover, $H(\sigma )$ is an $H$-bicomodule algebra via
$\Delta $. \par \smallskip

It is well-known that in general the set $Z^2 (H)$ of left
2-cocycles is not closed under convolution. One of the main
features of lazy 2-cocycles is that the set $Z^2_L (H)$ is closed
under convolution, and that the convolution inverse of an element
$\sigma \in Z^2_L (H)$ is again a lazy 2-cocycle, so $Z^2_L (H)$
is a group under convolution. In particular, a lazy 2-cocycle is
also a right 2-cocycle. Consider now the map $D^1:Reg^1
(H)\rightarrow Reg^2 (H)$, $D^1(\gamma )(h, h')=\gamma
(h_1)\gamma (h'_1)\gamma ^{-1}(h_2h'_2)$, for all $h, h'\in H$.
Then, by \cite{bc}, the map $D^1$ induces a group morphism
$Reg^1_L (H)\rightarrow Z^2_L (H)$, whose image is contained in
the centre of $Z^2_L (H)$; denote by $B^2_L (H)$ this central
subgroup $D^1(Reg^1_L (H))$ of $Z^2_L (H)$ (its elements are
called {\it lazy 2-coboundaries}). Finally, define the {\it second
lazy cohomology group} $H^2_L (H)=Z^2_L (H)/B^2_L (H)$ (most
likely nonabelian in general). Lazy 2-cocycles belonging to the
same class in $H^2_L (H)$ (we call them {\it lazy cohomologous})
are in particular cohomologous in the sense recalled before.

\section{Some properties of lazy 2-cocycles}\label{sec2}
\setcounter{equation}{0}

The aim of this section is to give some general properties of
lazy 2-cocycles needed in the next sections although they could
also be of independent interest. \vspace{2mm}

Let $\sigma :H\ot H\rightarrow k$ be a normalized and convolution
invertible left 2-cocycle. It is well-known (see \cite{montg},
\cite{c}) that the following formulae hold:
\begin{eqnarray}
&&\sigma (h_1, S(h_2))\sigma ^{-1}(S(h_3), h_4)=\varepsilon (h),
\label{n+1}\\
&&\sigma (S^{-1}(h_2), h_1)\sigma ^{-1}(h_4,
S^{-1}(h_3))=\varepsilon (h), \label{n+2}
\end{eqnarray}
for all $h\in H$, but in general we see no reason to have formulae of the type
\begin{eqnarray}
&&\sigma (h_1, S(h_2))=\sigma (S(h_1), h_2), \label{n+3}\\
&&\sigma (S^{-1}(h_2), h_1)=\sigma (h_2, S^{-1}(h_1)), \label{n+4}
\end{eqnarray}
even if for group algebras these formulae are true and
well-known. We have searched through the literature to find an
explicit counterexample, but we could not find any, so we are
going to provide here one. The Hopf algebra will be the Taft Hopf
algebra $H_9$ of dimension 9. \vspace{2mm}

Recall that $H_9=k\langle X, Y \vert
X^3=1,\;Y^3=0,\;YX=qXY\rangle$, where $q$ is a primitive 3-rd
root of unity, $\Delta (X)=X\ot X$, $\Delta (Y)=1\ot Y+Y\ot X$,
$S(X)=X^2$, $S(Y)=-q^2X^2Y$. The cleft extensions for any
$H_{n^2}$ have been classified in \cite{mas}, \cite{dt}; we use
here the form in \cite{dt}. We will construct a certain
$H_9$-cleft datum over $k$ (in the terminology of \cite{dt}).
Namely, in the notation of \cite[Theorem 3.5]{dt}, we choose
$F=id_k$, $D=0$ and $\alpha $, $\beta $, $\gamma $$\;\in k$ with
$\alpha , \gamma \neq 0$. Then, also in the notation of
\cite{dt}, one computes easily that:
\begin{eqnarray*}
&&\gamma _2=(1+q)\gamma ,\\
&&\gamma _3=(1+q+q^2)\gamma =0,\\
&&\theta_2=(1+q)\gamma^2, \\
&&\theta _3=D(\theta _2)+\theta _2\gamma _3=0, \\
&&D_3=0.
\end{eqnarray*}
Using these formulae, one can see that the conditions (1)-(9) in
\cite[Theorem 3.5]{dt}, are satisfied, so indeed $(id, 0, \alpha
, \beta, \gamma )$ is an $H_9$-cleft datum. The table for the
left 2-cocycle corresponding to any $H_9$-cleft datum is given in
\cite[Example 3.6]{dt}. For our datum, we get from the table:
\begin{eqnarray*}
&&\sigma (Y, XY)=0, \;\;\sigma (Y^2, X)=\theta _2\alpha =\gamma
^2(1+q)\alpha ,\\
&&\sigma (X^2Y, XY)=0, \;\;\sigma (XY^2, X^2)=FD(\gamma _2\alpha
)=0.
\end{eqnarray*}
From the equalities
\begin{eqnarray*}
&&\Delta (Y^2)=1\ot Y^2+(1+q)Y\ot XY+Y^2\ot X^2,\\
&&S(Y^2)=XY^2, \;S(XY)=-qXY, \;S(X^2)=X,
\end{eqnarray*}
we compute $\sigma (h_1, S(h_2))$ and $\sigma (S(h_1), h_2)$ for
the element $h:=Y^2$ and we obtain:
\begin{eqnarray*}
&&\sigma (h_1, S(h_2))=(1+q)\sigma (Y, -qXY)+\sigma (Y^2,
X)=\gamma ^2(1+q)\alpha \neq 0, \\
&&\sigma (S(h_1), h_2)=(1+q)\sigma (-q^2X^2Y, XY)+\sigma (XY^2,
X^2)=0,
\end{eqnarray*}
so the two terms cannot be equal.\vspace{2mm}

However, we have the following very useful result.

\begin{lemma}
If $\sigma $ is lazy, then formulae (\ref{n+3}) and (\ref{n+4})
hold.
\end{lemma}
\begin{proof}
Since $\sigma $ is lazy, the left cocycle condition can be
written as
\begin{eqnarray*}
&&\sigma (a_1, b_1)\sigma (a_2b_2, c)=\sigma (a, b_1c_1) \sigma
(b_2, c_2).
\end{eqnarray*}
By taking $a=h_1$, $b=S(h_2)$, $c=h_3$ in this formula, we obtain
(\ref{n+3}). Since $\sigma $ is lazy, it is also a right
2-cocycle, and the right cocycle condition can be written, using
the laziness of $\sigma$, as
\begin{eqnarray*}
&&\sigma (a_1b_1, c)\sigma (a_2, b_2)=\sigma (b_1, c_1)\sigma (a,
b_2c_2).
\end{eqnarray*}
By taking in this formula $a=h_3$, $b=S^{-1}(h_2)$, $c=h_1$, we
obtain (\ref{n+4}).
\end{proof}

We give now some more useful formulae.

\begin{lemma}
Let $\sigma :H\ot H\rightarrow k$ be a normalized and convolution
invertible lazy 2-cocycle. Then we have:
\begin{eqnarray}
&&\sigma ^{-1}(h_3, S^{-1}(h_2))h_4S^{-1}(h_1)=\sigma ^{-1}(h_2,
S^{-1}(h_1))1, \label{n+5}\\
&&\sigma ^{-1}(S^{-1}(h_3),h_2)S^{-1}(h_4)h_1=\sigma
^{-1}(S^{-1}(h_2), h_1)1, \label{n+10}\\
&&\sigma ^{-1}(S(h_2),h_3)S(h_1)h_4=\sigma
^{-1}(S(h_1), h_2)1, \label{n+6}\\
&&\sigma ^{-1}(S(h_2),h_3)S(h_1)=\sigma
^{-1}(S(h_1), h_2)S(h_3), \label{n+7}\\
&&\sigma ^{-1}(h_2, S(h_3))h_1S(h_4)=\sigma ^{-1}(h_1,
S(h_2))1, \label{n+8}\\
&&\sigma ^{-1}(S(h_2),h_3)h_4S^{-1}(h_1)=\sigma ^{-1}(S(h_1),
h_2)1. \label{n+9}
\end{eqnarray}
\end{lemma}
\begin{proof}
For (\ref{n+5}), apply the lazy condition to the elements $h_2$
and $S^{-1}(h_1)$; for (\ref{n+10}), apply the lazy condition to
the elements $S^{-1}(h_2)$ and $h_1$; for (\ref{n+6}), apply the
lazy condition to the elements $S(h_1)$ and $h_2$; (\ref{n+7}) is
obtained from (\ref{n+6}) by making convolution to the right with
$S$; for (\ref{n+8}), apply the lazy condition to the elements
$h_1$ and $S(h_2)$; finally, (\ref{n+9}) is obtained from
(\ref{n+8}) by using (\ref{n+3}) and then applying $S^{-1}$.
\end{proof}

Let $\sigma :H\ot H\rightarrow k$ be a normalized and convolution
invertible left 2-cocycle. Let us recall from \cite{bc} that the
linear map $\phi_{\sigma }:$ $_{\sigma}H\rightarrow H_{\sigma
^{-1}}$ defined by
\begin{eqnarray}
&&\phi _{\sigma}(h)=\sigma (h_1, S(h_2))S(h_3) \label{julien}
\end{eqnarray}
is an algebra antimorphism, and moreover it satisfies, for all
$h\in H$: $$\phi_{\sigma }(h_1)\cdot _{\sigma
^{-1}}h_2=\varepsilon (h)1= h_1\cdot _{\sigma ^{-1}}\phi_{\sigma
}(h_2).$$ Also, let us recall from \cite{c} the maps $S_1,
S_2:H\rightarrow H$ given for all $h\in H$ by
\begin{eqnarray}
&&S_1(h)=\sigma ^{-1}(S(h_2), h_3)S(h_1), \label{4} \\
&&S_2(h)=\sigma ^{-1}(h_3, S^{-1}(h_2))S^{-1}(h_1). \label{5}
\end{eqnarray}
From (\ref{n+3}) and (\ref{n+7}) it follows immediately that:

\begin{proposition}
If $\sigma $ is lazy, then $S_1=\phi_{\sigma ^{-1}}$.
\end{proposition}

There exists also a relation between $S_2$ and $\phi_{\sigma }$,
which holds in general.

\begin{proposition}
If $\sigma $ is a normalized and convolution invertible left
2-cocycle on $H$, then $S_2$ is the composition inverse of
$\phi_{\sigma }$. In particular, it follows that $\phi_{\sigma }$
is bijective.
\end{proposition}
\begin{proof}
That $S_2\circ \phi_{\sigma }=id$ and $\phi_{\sigma }\circ S_2=id$
reduce respectively to formulae (\ref{n+1}) and (\ref{n+2}).
\end{proof}

If $\sigma $ is lazy, since $S_1=\phi_{\sigma ^{-1}}$ and $S_2$
is the composition inverse of $\phi_{\sigma }$, from the
properties of $\phi_{\sigma }$ we obtain:

\begin{proposition}
If $\sigma $ is lazy, then
$S_1, S_2:H(\sigma ^{-1})\rightarrow H(\sigma )$ are algebra
antiisomorphisms, and we have, for all $h\in H$:
\begin{eqnarray}
&&S_1(h_1)\cdot _{\sigma }h_2=\varepsilon (h)1=h_1\cdot _{\sigma
}S_1(h_2), \label{6}\\
&&S_2(h_2)\cdot _{\sigma }h_1=\varepsilon (h)1=h_2\cdot _{\sigma
}S_2(h_1). \label{7}
\end{eqnarray}
\end{proposition}

Let us note that (\ref{6}) and (\ref{7}) appear also in \cite{c},
in a slightly different form, and they actually hold for any left
2-cocycle, not necessarily lazy.

\begin{proposition}
Let $\sigma $ be a normalized and convolution invertible left
2-cocycle on $H$. Then we have, for all $h\in H$:
\begin{eqnarray}
&&\Delta (S_1(h))=S_1(h_2)\ot S(h_1), \label{8}\\
&&\Delta (S_2(h))=S_2(h_2)\ot S^{-1}(h_1). \label{9}
\end{eqnarray}
\end{proposition}
\begin{proof}
An easy computation.
\end{proof}

\section{Extending lazy 2-cocycles to a Drinfeld double}\label{sec3}
\setcounter{equation}{0}

Throughout this section, $H$ will be a finite dimensional Hopf
algebra and we will denote the Drinfeld double of $H$ by $D(H)$.
A complete description of $H^2_L(D(H))$ in terms of $H^2_L(H)$
and $H^2_L(H^*)$ was given in \cite{bc}. In particular, it
follows from \cite{bc} that if $\sigma $ is a normalized and
convolution invertible lazy 2-cocycle on $H$, then it can be
extended to a normalized and convolution invertible lazy
2-cocycle $\overline{\sigma}$ on $D(H)$. In this section we provide an
alternative approach to the problem of extending a lazy 2-cocycle from
$H$ to $D(H)$, based on the so-called
diagonal crossed product construction. The results in this section will be
also used in Section \ref{sec5}.
\vspace{2mm}

Recall that the Drinfeld double of $H$ is a quasitriangular Hopf
algebra realized on the $k$-linear space $H^*\ot H$; its
coalgebra structure is $H^{*cop}\ot H$ and the algebra
structure is given by $$(p\ot h)(q\ot l)=p(h_1\rightharpoonup
q\leftharpoonup S^{-1}(h_3))\ot h_2l,$$ for all $p, q\in H^*$ and
$h, l\in H$, where $\rightharpoonup $ and $\leftharpoonup $ are
the left and right regular actions of $H$ on $H^*$ given by
$(h\rightharpoonup p)(l)=p(lh)$ and $(p\leftharpoonup
h)(l)=p(hl)$ for all $h, l\in H$ and $p\in H^*$. Let now $A$ be
an $H$-bicomodule algebra, with comodule structures $A\rightarrow
A\ot H$, $a\mapsto a_{<0>}\ot a_{<1>}$ and $A\rightarrow H\ot A$,
$a\mapsto a_{[-1]}\ot a_{[0]}$, and denote, for $a\in A$,
$$a_{\{-1\}}\ot a_{\{0\}}\ot a_{\{1\}}=a_{<0>_{[-1]}}\ot
a_{<0>_{[0]}}\ot a_{<1>}=a_{[-1]}\ot a_{[0]_{<0>}}\ot
a_{[0]_{<1>}},$$ as an element in $H\ot A\ot H$. Recall from
\cite{hn1} that the (left) {\it diagonal crossed product}
$H^*\bowtie A$ is equal to $H^*\ot A$ as a $k$-space, but with
multiplication given by: $$(p\bowtie a)(q\bowtie
b)=p(a_{\{-1\}}\rightharpoonup q\leftharpoonup
S^{-1}(a_{\{1\}}))\bowtie a_{\{0\}}b,$$ for all $a, b\in A$ and
$p, q\in H^*$, and with unit $\varepsilon _H\bowtie 1_A$. The
space $H^*\bowtie A$ becomes a $D(H)$-bicomodule algebra, with
structures
\begin{eqnarray*}
&&H^*\bowtie A\rightarrow (H^*\bowtie A)\ot D(H), \ p\bowtie a\mapsto
(p_2\bowtie a_{<0>})\ot (p_1\ot a_{<1>}), \\
&&H^*\bowtie A\rightarrow D(H)\ot (H^*\bowtie A), \ p\bowtie
a\mapsto (p_2\ot a_{[-1]})\ot (p_1\bowtie a_{[0]}),
\end{eqnarray*}
for all $p\in H^*$, $a\in A$. If $A=H$ then $H^*\bowtie A$ is just
$D(H)$, with bicomodule algebra structure over itself given by
its comultiplication. It is well-known (see \cite{dt2}) that the
Drinfeld double can be expressed as a twisting of $H^{*cop}\ot
H$. Similarly, using the framework and notation of \cite{sch1},
one can prove that
$$\begin{array}{ccc}
                      & \hspace{-5pt} \tau & \vspace{-3pt}\\
H^*\bowtie A=H^{*cop} & \hspace{-5pt} \# & \hspace{-5pt} A  \vspace{-3pt}\\
 & \hspace{-5pt} \tau &
\end{array}$$
where $\tau :H\ot H^{*cop}\rightarrow k$ is the skew-pairing
given by $\tau (h, p)=p(h)$. \\

Let $\sigma:H \otimes H \rightarrow k$ be a normalized and
invertible lazy $2$-cocyle. Either as a consequence of the proof
in \cite{bc}, or by direct means, one can see that the extended
cocycle $\overline{\sigma}:D(H) \otimes D(H) \rightarrow k$ and
its convolution inverse are given by the formulae
\begin{eqnarray}
&&\overline{\sigma }\;(p\ot h, q\ot l)=p(1)q(S^{-1}(h_3)h_1)\sigma
(h_2, l),  \label{1.4}\\
&&\overline{\sigma }^{\;-1}(\;p\ot h, q\ot
l)=p(1)q(S^{-1}(h_3)h_1)\sigma ^{-1}(h_2, l),
\end{eqnarray}
for all $p, q\in H^*$ and $h, l\in H$. \\

In view of the above description of the diagonal crossed product
as a twisting and of the nature of the proof for the description
of $H^2_L(D(H))$ in \cite{bc}, it is likely that the following
result can be proved using the approach in \cite{bc}. But we
prefer to give a direct proof, because this is how we discovered
it (actually, how we got the formula (\ref{1.4}) for
$\overline{\sigma }$).
\begin{proposition}\label{exte}
Let $\sigma :H\ot H\rightarrow k$ be a normalized and convolution
invertible lazy 2-cocycle. Consider the $H$-bicomodule algebra
$H(\sigma )$. Then $H^*\bowtie H(\sigma )=D(H)(\overline{\sigma
}),$ as $D(H)$-bicomodule algebras. Moreover, $\overline{\sigma
}$ is unique with this property.
\end{proposition}
\begin{proof}
We compute the multiplications in the two algebras and show that
they coincide.
$$\begin{array}{l}
(p\bowtie h)(q\bowtie l)=p(h_1\rightharpoonup q\leftharpoonup
S^{-1}(h_3))\bowtie h_2\cdot _{\sigma }l \vspace{3pt} \\
\hspace{1.5cm} =\sigma (h_2, l_1)p(h_1\rightharpoonup
q\leftharpoonup S^{-1}(h_4))\bowtie h_3l_2, \vspace{3pt} \\
\hspace{1.5cm} =p(h_4S^{-1}(h_3)h_1\rightharpoonup q\leftharpoonup
S^{-1}(h_6))\sigma (h_2, l_1)\ot h_5l_2 \vspace{3pt} \\
\hspace{1.5cm} =p(h_4\rightharpoonup q_1\leftharpoonup
S^{-1}(h_6))q_2(S^{-1}(h_3)h_1)\sigma (h_2, l_1)\ot h_5l_2 \vspace{3pt} \\
\hspace{1.5cm} =p_2(1)q_2(S^{-1}(h_{(1,3)})h_{(1, 1)})\sigma
(h_{(1, 2)}, l_1)p_1(h_{(2, 1)}\rightharpoonup q_1\leftharpoonup
S^{-1}(h_{(2,3)}))\ot h_{(2, 2)}l_2 \vspace{3pt} \\
\hspace{1.5cm} =\overline{\sigma }(p_2\ot h_1, q_2\ot l_1)(p_1\ot
h_2)(q_1\ot l_2) \vspace{3pt} \\
\hspace{1.5cm} =(p\ot h)\cdot _{\overline{\sigma }}(q\ot l).
\end{array}$$
Clearly $H^*\bowtie H(\sigma)$ and $D(H)(\overline{\sigma })$ have
the same $D(H)$-bicomodule structure. For the uniqueness of
$\overline{\sigma }$, we write down the fact that the
multiplications in $H^*\bowtie H(\sigma )$ and
$D(H)(\overline{\sigma })$ coincide, then we evaluate this
equality on $1\ot \varepsilon $ and we obtain that
$\overline{\sigma }$ has to be given by (\ref{1.4}).
\end{proof}

It was proved in \cite{bc} that $H^2_L(H)$ can be embedded as a
subgroup in $Bigal (H)$, the group of Bigalois objects of $H$
introduced in \cite{sch2}, \cite{voz1}.
\begin{proposition}
The map $A\mapsto H^*\bowtie A$ gives an embedding of $Bigal (H)$
into $Bigal (D(H))$, whose restriction to $H^2_L(H)$ is the
embedding of $H^2_L(H)$ into $H^2_L(D(H))$ from \cite{bc}.
\end{proposition}
\begin{proof}
The fact that the map $A\mapsto H^*\bowtie A$ gives the desired
embedding between Bigalois groups is contained, even if not
explicitly stated, in Schauenburg's paper \cite{sch1}, and the
compatibility between the two embeddings, at the levels of
Bigalois groups and lazy cohomologies, follows from the
compatibility between the proof in \cite{sch1} and the one in
\cite{bc}.
\end{proof}

The antipode of $D(H)$ is given by the formula
\begin{eqnarray*}
&&S_{D(H)}(p\ot h)=(\varepsilon \ot S(h))(S^{*-1}(p)\ot 1),
\end{eqnarray*}
for all $h\in H$, $p\in H^*$. Denote $S_{D(H)}$ by
$\overline{S}$. One can easily check that its inverse is given by
\begin{eqnarray*}
&&S_{D(H)}^{-1}(p\ot h)=(\varepsilon \ot S^{-1}(h))(S^*(p)\ot 1).
\end{eqnarray*}

Let now $\sigma :H\ot H\rightarrow k$ be a normalized and
convolution invertible lazy 2-cocycle, and $\overline{\sigma }$
its extension to $D(H)$, given by the formula (\ref{1.4}). Denote
by $S_1, S_2: H\rightarrow H$ the maps given by the formulae
(\ref{4}), (\ref{5}), and by $\overline{S}_1, \overline{S}_2:
D(H)\rightarrow D(H)$ the analogous maps for $D(H)$ corresponding
to $\overline{\sigma }$, that is:
\begin{eqnarray}
&&\overline{S}_1(p\ot h)=\overline{\sigma
}\;^{-1}(\overline{S}((p\ot h)_2), (p\ot h)_3)\overline{S}((p\ot
h)_1), \label{sb1}\\
&&\overline{S}_2(p\ot h)=\overline{\sigma }\;^{-1}((p\ot h)_3,
\overline{S}\;^{-1}((p\ot h)_2))\overline{S}\;^{-1}((p\ot h)_1).
\label{sb2}
\end{eqnarray}

The following result will be needed in a subsequent section.

\begin{proposition}
$\overline{S}_1$ and $\overline{S}_2$ can be computed as:
\begin{eqnarray}
&&\overline{S}_1(p\ot h)=(\varepsilon \ot S_1(h))(S^{*-1}(p)\ot
1),\label{sd1}\\
&&\overline{S}_2(p\ot h)=(\varepsilon \ot S_2(h))(S^*(p)\ot 1),
\label{sd2}
\end{eqnarray}
for all $h\in H$, $p\in H^*$.
\end{proposition}
\begin{proof}
We give the proof for $\overline{S}_1$, the one for
$\overline{S}_2$ is similar (but for $\overline{S}_2$ one has to
use the formula (\ref{n+5})). We compute:
$$\begin{array}{ll}
\overline{S}_1(p\ot h) & =\overline{\sigma }\;^{-1}(\overline{S}
(p_2\ot h_2), p_1\ot h_3)\overline{S} (p_3\ot h_1) \vspace{3pt} \\
 & =\overline{\sigma }\;^{-1}((\varepsilon \ot
S(h_2))(S^{*-1} (p_2)\ot 1), p_1\ot h_3)\overline{S} (p_3\ot h_1) \vspace{3pt} \\
 & =\overline{\sigma
}\;^{-1}(S(h_2)_1\rightharpoonup S^{*-1}(p_2)\leftharpoonup
S^{-1}(S(h_2)_3)\ot S(h_2)_2, p_1\ot
h_3)\overline{S} (p_3\ot h_1) \vspace{3pt} \\
 & =\overline{\sigma }\;^{-1}(S(h_4)\rightharpoonup
S^{*-1}(p_2)\leftharpoonup h_2\ot S(h_3), p_1\ot h_5)\overline{S}
(p_3\ot h_1) \vspace{3pt} \\
 & =(S(h_4)\rightharpoonup
S^{*-1}(p_2)\leftharpoonup
h_2)(1)p_1(S^{-1}(S(h_3)_3)S(h_3)_1)\sigma ^{-1}(S(h_3)_2, h_5)
\overline{S}(p_3\ot h_1) \vspace{3pt} \\
 & =S^{*-1}(p_2)(h_2S(h_6))p_1(h_3S(h_5))
\sigma^{-1}(S(h_4), h_7)\overline{S}(p_3\ot h_1) \vspace{3pt}\\
 & =p_1(h_3S(h_5)h_6S^{-1}(h_2))\sigma ^{-1}(S(h_4),
h_7)\overline{S}(p_2\ot h_1) \vspace{3pt} \\
 & =\sigma ^{-1}(S(h_2), h_3)(\varepsilon \ot
S(h_1))(S^{*-1}(p)\ot 1) \vspace{3pt}\\
 & =(\varepsilon \ot S_1(h))(S^{*-1}(p)\ot 1),
\end{array}$$
which was what we had to prove.
\end{proof}
\begin{remark}{\em
Using either the formula for $\overline{\sigma }$ or the identification
$H^*\bowtie H(\sigma )=D(H)(\overline{\sigma })$, one can easily check
that
\begin{eqnarray}
&&(\varepsilon \ot h)\cdot _{\overline{\sigma }}(p\ot
1)=(\varepsilon \ot h)(p\ot 1), \label{lulu}
\end{eqnarray}
for all $h\in H$ and $p\in H^*$ (we will use this later).}
\end{remark}

\section{Extending (lazy) 2-cocycles to a Radford biproduct}\label{sec4}
\setcounter{equation}{0}

For a Hopf algebra $H$ and a Hopf algebra $B$ in the category of
left Yetter-Drinfeld modules $_H^H{\cal YD}$ it is possible to
construct the Radford biproduct Hopf algebra $B \times H$. A
second lazy cohomology group $H^2_L(B)$ can be defined for $B$
inside the category $_H^H{\cal YD}$. In this section we find out
a relation between $H^2_L(B)$ and $H^2_L(B \times H)$.
\vspace{2mm}

We start by recalling from \cite{rad} the construction of a
Radford biproduct. Let $H$ be a bialgebra and $B$ a vector space
such that $(B, 1_B)$ is an algebra (with multiplication denoted by
$b\ot c\mapsto bc$ for all $b, c\in B$) and $(B, \Delta _B,
\varepsilon _B)$ is a coalgebra. The pair $(H, B)$ is called {\it
admissible} if $B$ is endowed with a left $H$-module structure
(denoted by $h\ot b\mapsto h\cdot b)$ and with a left
$H$-comodule structure (denoted by $b\mapsto b^{(-1)}\ot
b^{(0)}\in H\ot B)$ such that:\\

(1) $B$ is a left $H$-module algebra;

(2) $B$ is a left $H$-comodule algebra;

(3) $B$ is a left $H$-comodule coalgebra, that is, for all $b\in
B$:
\begin{eqnarray}
&&b_1^{(-1)}b_2^{(-1)}\ot b_1^{(0)}\ot b_2^{(0)}=b^{(-1)}\ot (b^{(0)})_1
\ot (b^{(0)})_2, \label{r1} \\
&&b^{(-1)}\varepsilon _B(b^{(0)})=\varepsilon _B(b)1_H. \label{r2}
\end{eqnarray}

(4) $B$ is a left $H$-module coalgebra, that is, for all $h\in H$
and $b\in B$:
\begin{eqnarray}
&&\Delta _B(h\cdot b)=h_1\cdot b_1\ot h_2\cdot b_2, \label{r3} \\
&&\varepsilon _B(h\cdot b)=\varepsilon _H(h)\varepsilon _B(b). \label{r4}
\end{eqnarray}

(5) $\varepsilon _B$ is an algebra map and $\Delta _B(1_B)=1_B\ot
1_B$;

(6) The following relations hold for all $h\in H$ and $b, c\in B$:
\begin{eqnarray}
&&\Delta _B(bc)=b_1(b_2^{(-1)}\cdot c_1)\ot b_2^{(0)}c_2, \label{r5} \\
&&(h_1\cdot b)^{(-1)}h_2\ot (h_1\cdot b)^{(0)}=h_1b^{(-1)}\ot
h_2\cdot b^{(0)}. \label{r6}
\end{eqnarray}

If $(H, B)$ is an admissible pair, then we know from \cite{rad}
that the smash product algebra structure and smash coproduct
coalgebra structure on $B\ot H$ afford $B\ot H$ a bialgebra
structure, denoted by $B\times H$ and called the {\it smash
biproduct} or {\it Radford biproduct}. Its comultiplication is
given by
\begin{eqnarray}
&&\Delta (b\times h)=(b_1\times b_2^{(-1)}h_1)\ot (b_2^{(0)}\times h_2),
\label{r7}
\end{eqnarray}
for all $b\in B$, $h\in H$, and its counit is $\varepsilon _B\ot
\varepsilon _H$. Let us record the following formula:
\begin{eqnarray}
&&\Delta _B(b(h\cdot c))=b_1(b_2^{(-1)}h_1\cdot c_1)\ot b_2^{(0)}
(h_2\cdot c_2), \label{r8}
\end{eqnarray}
for all $h\in H$ and $b, c\in B$, which follows immediately from
(\ref{r5}) and (\ref{r3}). If $H$ is a Hopf algebra with antipode
$S_H$ and $(H, B)$ is an admissible pair such that there exists
$S_B\in Hom (B, B)$ a convolution inverse for $id_B$, then
$B\times H$ is a Hopf algebra with antipode
\begin{eqnarray}
S(b\times h)=(1\times S_H(b^{(-1)}h))(S_B(b^{(0)})\times 1), \label{r9}
\end{eqnarray}
for all $h\in H$, $b\in B$. In this case, we will say that $(H,
B)$ is a {\it Hopf admissible pair}. For a Hopf algebra $H$, it is
well-known (see for instance \cite{montg}, \cite{majid}) that
$(H, B)$ being an admissible pair (respectively Hopf admissible
pair) is equivalent to $B$ being a bialgebra (respectively Hopf
algebra) in the Yetter-Drinfeld category $_H^H{\cal
YD}$.\vspace{2mm}

Recall now from \cite{tak} the so-called {\it generalized smash
product}. If $H$ is a bialgebra, $B$ a left $H$-module algebra
(with action $h\ot b\mapsto h\cdot b)$ and $A$ a left $H$-comodule
algebra (with coaction $a\mapsto a_{(-1)}\ot a_{(0)}\in H\ot A)$,
then on $B\ot A$ we have an associative algebra structure,
denoted by $B\gsm A$, with unit $1_B\gsm 1_A$ and multiplication
\begin{eqnarray}
&&(b\gsm a)(b'\gsm a')=b(a_{(-1)}\cdot b')\gsm a_{(0)}a', \label{r10}
\end{eqnarray}
for all $b, b'\in B$ and $a, a'\in A$. \par \vspace{2mm}

As we have seen before, the relation between a Drinfeld double
and a diagonal crossed product is that the diagonal crossed
product becomes a bicomodule algebra over the Drinfeld double.
The next result shows that a similar relation exists between a
Radford biproduct and a generalized smash product.

\begin{proposition}
If $(H, B)$ is an admissible pair and $A$ is a left $H$-comodule algebra,
then $B\gsm A$ becomes a left $B\times H$-comodule algebra, with coaction
$$\lambda :B\gsm A\rightarrow (B\times H)\ot (B\gsm A), \
\lambda (b\gsm a)=(b_1\times b_2^{(-1)}a_{(-1)})\ot (b_2^{(0)}
\gsm a_{(0)}),$$ for all $b\in B$ and $a\in A$.
\end{proposition}
\begin{proof}
We prove first that $(B\gsm A, \lambda )$ is a left $B\times
H$-comodule (for this part we only need $A$ to be a left
$H$-comodule). We compute:
$$\begin{array}{l}
(id \ot \lambda )(\lambda (b\gsm a))= (b_1\times
b_2^{(-1)}a_{(-1)})\ot ((b_2^{(0)})_1\times (b_2^{(0)})_2^{(-1)}
a_{(0)_{(-1)}})\ot ((b_2^{(0)})_2^{(0)}\gsm a_{(0)_{(0)}}) \vspace{3pt} \\
\hspace{1.5cm}=(b_1\times b_2^{(-1)}a_{(-1)_1})\ot
((b_2^{(0)})_1\times (b_2^{(0)})_2^{(-1)}a_{(-1)_2})\ot ((b_2^{(0)})_2^{(0)}
\gsm a_{(0)})  \vspace{3pt} \\
\hspace{1.5cm}=(b_1\times b_2^{(-1)}b_3^{(-1)}a_{(-1)_1})\ot
(b_2^{(0)}\times b_3^{(0)^{(-1)}}a_{(-1)_2})\ot
(b_3^{(0)^{(0)}}\gsm a_{(0)})\hfill (by \;\; (\ref{r1}))  \vspace{3pt} \\
\hspace{1.5cm}=(b_1\times b_2^{(-1)}(b_3^{(-1)})_1a_{(-1)_1})\ot
(b_2^{(0)}\times (b_3^{(-1)})_2a_{(-1)_2})\ot (b_3^{(0)}\gsm a_{(0)})  \vspace{3pt} \\
\hspace{1.5cm}=(b_{(1, 1)}\times
b_{(1,2)}^{(-1)}(b_2^{(-1)})_1a_{(-1)_1})\ot (b_{(1,
2)}^{(0)}\times (b_2^{(-1)})_2a_{(-1)_2})\ot (b_2^{(0)}\gsm a_{(0)}) \vspace{3pt} \\
\hspace{1.5cm}=\Delta (b_1\times b_2^{(-1)}a_{(-1)})\ot (b_2^{(0)}\gsm a_{(0)}) \vspace{3pt} \\
\hspace{1.5cm}=(\Delta \ot id)(\lambda (b \gsm a)).
\end{array}$$
Then obviously we have that $(\varepsilon \ot id)\lambda =id$, so
$B\gsm A$ is indeed a left $B\times H$-comodule. We proceed to
show that $\lambda $ is an algebra map. First, by (5), we have
$\lambda (1\gsm 1)=(1\times 1)\ot (1\gsm 1)$. For $b, b'\in B$
and $a, a'\in A$ we have:
$$\begin{array}{l}
\lambda ((b\gsm a)(b'\gsm a'))=\lambda (b(a_{(-1)}\cdot b')\gsm a_{(0)}a') \vspace{3pt} \\
\hspace{1.2cm}=((b(a_{(-1)}\cdot b'))_1\times (b(a_{(-1)}\cdot
b'))_2^{(-1)}(a_{(0)}a')_{(-1)})\ot ((b(a_{(-1)}\cdot b'))_2^{(0)}\gsm (a_{(0)}a')_{(0)}) \vspace{3pt}\\
\hspace{1.2cm} =(b_1(b_2^{(-1)}a_{(-1)_1}\cdot b'_1)\times
b_2^{(0)^{(-1)}}(a_{(-1)_2}\cdot b'_2)^{(-1)}a_{(0)_{(-1)}}a'_{(-1)})\ot \vspace{3pt} \\
\hspace{1.2cm} \;\;\;\;\;\;\;\;\ot
\;(b_2^{(0)^{(0)}}(a_{(-1)_2}\cdot b'_2)^{(0)}\gsm
a_{(0)_{(0)}}a'_{(0)}) \hfill (by \;\; (\ref{r8})) \vspace{3pt} \\
\hspace{1.2cm}=(b_1(b_2^{(-1)}a_{(-1)_1}\cdot b'_1)\times
b_2^{(0)^{(-1)}}(a_{(-1)_2}\cdot b'_2)^{(-1)}a_{(-1)_3}a'_{(-1)})\ot \vspace{3pt} \\
\hspace{1.2cm}\;\;\;\;\;\;\;\;\ot \;(b_2^{(0)^{(0)}}(a_{(-1)_2}\cdot b'_2)^{(0)}\gsm a_{(0)}a'_{(0)}) \vspace{3pt} \\
\hspace{1.2cm}=(b_1(b_2^{(-1)}a_{(-1)_1}\cdot b'_1)\times
b_2^{(0)^{(-1)}}a_{(-1)_2}b_2^{'(-1)}a'_{(-1)})\ot \vspace{3pt} \\
\hspace{1.2cm}\;\;\;\;\;\;\;\;\;\ot \;(b_2^{(0)^{(0)}}
(a_{(-1)_3}\cdot b_2^{'(0)})\gsm a_{(0)}a'_{(0)}) \hfill
(by \;\;(\ref{r6})) \vspace{3pt} \\
\hspace{1.2cm}=(b_1(b_2^{(-1)}a_{(-1)_1}\cdot b'_1)\times
b_2^{(0)^{(-1)}}a_{(-1)_2} b_2^{'(-1)}a'_{(-1)})\ot
(b_2^{(0)^{(0)}}(a_{(0)_{(-1)}}\cdot b_2^{'(0)})\gsm a_{(0)_{(0)}}a'_{(0)}) \vspace{3pt} \\
\hspace{1.2cm}=(b_1((b_2^{(-1)})_1a_{(-1)_1}\cdot b'_1)\times
(b_2^{(-1)})_2a_{(-1)_2} b_2^{'(-1)}a'_{(-1)})\ot (b_2^{(0)}
(a_{(0)_{(-1)}}\cdot b_2^{'(0)})\gsm a_{(0)_{(0)}}a'_{(0)}) \vspace{3pt}\\
\hspace{1.2cm}=\lambda (b\gsm a)\lambda (b'\gsm a'),
\end{array}$$
and the proof is finished.
\end{proof}

Now, let $(H, B)$ be an admissible pair and $\sigma :H\ot
H\rightarrow k$ a normalized and convolution invertible right
2-cocycle, so that we can consider $H_{\sigma }$, which is a left
$H$-comodule algebra, and we can make $B\gsm H_{\sigma }$, which,
by the above proposition, becomes a left $B\times H$-comodule
algebra.

\begin{proposition}
With notation as above, the map $\tilde{\sigma }:(B\times H)\ot
(B\times H)\rightarrow k$ defined by
\begin{eqnarray*}
&&\tilde{\sigma }(b\times h, b'\times h')=\varepsilon
_B(b)\varepsilon _B(b') \sigma (h, h'), \;\;\;\;\;\forall \;\;b,
b'\in B \;\;and\;\;h, h'\in H,
\end{eqnarray*}
is a normalized and convolution invertible right 2-cocycle on $B\times H$,
and we have $(B\times H)_{\tilde {\sigma }}=B\gsm H_{\sigma }$ as left
$B\times H$-comodule algebras. Moreover, $\tilde{\sigma }$ is unique with
this property.
\end{proposition}
\begin{proof}
We have:
$$\begin{array}{lll}
(b\gsm h)(b'\gsm h') & = b(h_1\cdot b')\gsm h_2\cdot _{\sigma }h' & \vspace{3pt} \\
& =b(h_1\cdot b')\gsm h_2h'_1 \sigma (h_3, h'_2) & \vspace{3pt} \\
& =(b\times h_1)(b'\times h'_1)\sigma (h_2, h'_2) &  \vspace{3pt} \\
& = (b_1\times b_2^{(-1)}h_1)(b'_1\times b_2^{'(-1)}h'_1)
\varepsilon _B(b_2^{(0)})\varepsilon _B(b_2^{'(0)})
\sigma (h_2, h'_2) & \quad (by \;\;(\ref{r2})) \vspace{3pt} \\
& = b_1\times b_2^{(-1)}h_1)(b'_1\times b_2^{'(-1)}h'_1)
\tilde{\sigma }(b_2^{(0)}\times h_2, b_2^{'(0)}\times h'_2) & \vspace{3pt} \\
& = (b\times h)_1(b'\times h')_1\tilde{\sigma }((b\times h)_2,
(b'\times h')_2). & \vspace{3pt} \\
& = (b\times h)\cdot _{\tilde{\sigma}} (b'\times h') &
\end{array}$$
So, the multiplication in $(B\times H)_{\tilde{\sigma }}$
coincides with the one in $B\gsm H_{\sigma }$, which is
associative, so $\tilde{\sigma }$ is automatically a right
2-cocycle, and we have $(B\times H)_{\tilde{\sigma }} =B\gsm
H_{\sigma }$ as algebras; it is obvious that they coincide also as
left $B\times H$-comodules, and is easy to prove that
$\tilde{\sigma }$ is normalized and convolution invertible. To
prove the uniqueness of $\tilde{\sigma }$, write that the
multiplications in $B\gsm H_{\sigma }$ and $(B\times
H)_{\tilde{\sigma }}$ coincide, apply $\varepsilon _B\ot
\varepsilon _H$ and get $\tilde{\sigma }(b\times h, b'\times
h')=\varepsilon _B(b)\varepsilon _B(b') \sigma (h, h')$.
\end{proof}

The map $\pi:B\times H \rightarrow H,\ b \times h \mapsto
\varepsilon(b) h$ is a Hopf algebra map. Observe that
$\tilde{\sigma}$ is just the cocycle obtained by pulling back
through the map $\pi$.

\begin{remark}{\em With notation as above, we have:
\begin{eqnarray*}
&&(b\times h)_1(b'\times h')_1\tilde{\sigma }((b\times h)_2,
(b'\times h')_2)=b(h_1\cdot b')\times h_2h'_1 \sigma (h_3, h'_2), \\
&&(b\times h)_2(b'\times h')_2\tilde{\sigma }((b\times h)_1,
(b'\times h')_1)=\sigma (b^{(-1)}h_1, b^{'(-1)}h'_1)b^{(0)}
(h_2\cdot b^{'(0)})\times h_3h'_2,
\end{eqnarray*}
for all $b, b'\in B$ and $h, h'\in H$. Assume that
$\tilde{\sigma}$ is lazy; then, by taking $b=b'=1$ above, we
obtain that $\sigma $ is lazy. Conversely, if $\sigma $ is lazy,
then $\tilde{\sigma }$ is lazy if and only if
\begin{eqnarray*}
&&\sigma (h_2, h')b(h_1\cdot b')=\sigma (b^{(-1)}h_1, b^{'(-1)}h')
b^{(0)}(h_2\cdot b^{'(0)}),
\end{eqnarray*}
for all $b, b'\in B$ and $h, h'\in H$, from which follow some necessary
conditions for the laziness of $\tilde{\sigma }$, such as
\begin{eqnarray*}
&&h\cdot b=\sigma (h_1, b^{(-1)})(h_2\cdot b^{(0)}),\\
&&bb'=\sigma (b^{(-1)}, b^{'(-1)})b^{(0)}b^{'(0)},
\end{eqnarray*}
for all $b, b'\in B$ and $h\in H$, which have no reason to hold in
general.}
\end{remark}

We study now the problem of extending (lazy) 2-cocycles from $B$
to $B\times H$.\vspace{2mm}

Let ${\cal C}$ be a braided monoidal category and $B$ a Hopf
algebra in ${\cal C}$. Then, just as if $B$ would be a usual Hopf
algebra, one can define 2-cocycles, crossed products, Galois
extensions, etc, for $B$ in ${\cal C}$, see for instance
\cite{zhang}, \cite{bespalov}. Also, one can define lazy
2-cocycles, lazy 2-coboundaries and the second lazy cohomology
group $H^2_L(B)=Z^2_L(B)/B^2_L(B)$. Here, we will only be
interested in the case when ${\cal C}$=$_H^H{\cal YD}$, the
category of left Yetter-Drinfeld modules over a Hopf algebra $H$,
and $B$ a Hopf algebra in $_H^H{\cal YD}$ (that is, $(H, B)$ is a
Hopf admissible pair, so $B\times H$ is a Hopf algebra). For this
category, one can prove by hand all the properties of lazy
2-cocycles that allow to define $H^2_L(B)$ (the most difficult is
to prove that the product of two lazy 2-cocycles is a left
2-cocycle--we will give an easy alternative proof of this fact at
the end of the section). \vspace{2mm}

If $M, N\in {}_H^H{\cal YD}$, then $M\ot N\in {}_H^H{\cal YD}$
with module structure $h\cdot (m\ot n)= h_1\cdot m\ot h_2\cdot n$
and comodule structure $m\ot n\mapsto m_{<-1>}n_{<-1>}\ot
(m_{<0>}\ot n_{<0>})$, where $m\mapsto m_{<-1>}\ot m_{<0>}$ and
$n\mapsto n_{<-1>}\ot n_{<0>}$ are the comodule structures of $M$
and $N$, and the braiding is given by
\begin{eqnarray}
&&c_{M, N}:M\ot N\rightarrow N\ot M, \;\;\;\;c_{M, N}(m\ot n)=m_{<-1>}
\cdot n\ot m_{<0>}. \label{b1}
\end{eqnarray}
Hence, the coalgebra structure of $B\ot B$ in $_H^H{\cal YD}$ is given by
\begin{eqnarray*}
\Delta _{B\ot B}(b\ot b')&=&(id \ot c_{B, B}\ot id)\circ (\Delta _B\ot
\Delta _B)(b\ot b')\\
&=&(b_1\ot b_2^{(-1)}\cdot b'_1)\ot (b_2^{(0)}\ot b'_2).
\end{eqnarray*}
So, if $\sigma ,\tau:B\ot B\rightarrow k$ are morphisms in
$_H^H{\cal YD}$, their convolution in $_H^H{\cal YD}$ is given by:
\begin{eqnarray}
&&(\sigma * \tau)(b\ot b')=\sigma (b_1\ot b_2^{(-1)}\cdot b'_1)
\tau (b_2^{(0)}\ot b'_2). \label{b3}
\end{eqnarray}

Let $\sigma :B\ot B\rightarrow k$ be a morphism in
$_H^H{\cal YD}$, that is, it satisfies the conditions:
\begin{eqnarray}
&&\sigma (h_1\cdot b\ot h_2\cdot b')=\varepsilon (h)\sigma (b\ot b'),
\label{b4} \\
&&\sigma (b^{(0)}\ot b^{'(0)})b^{(-1)}b^{'(-1)}=\sigma (b\ot b')1_H,
\label{b5}
\end{eqnarray}
for all $h\in H$ and $b, b'\in B$. Then $\sigma $ is a lazy element if it
satisfies the categorical laziness condition:
\begin{eqnarray}
&&\sigma (b_1\ot b_2^{(-1)}\cdot b'_1)b_2^{(0)}b'_2=\sigma (b_2^{(0)}\ot
b'_2)b_1(b_2^{(-1)}\cdot b'_1), \label{b6}
\end{eqnarray}
for all $b, b'\in B$.\vspace{2mm}

Let $\sigma :B\ot B\rightarrow k$ be a normalized left 2-cocycle in
$_H^H{\cal YD}$, that is $\sigma $ is a normalized morphism in
$_H^H{\cal YD}$ satisfying the categorical left 2-cocycle condition
\begin{eqnarray}
\sigma (a_1\ot a_2^{(-1)}\cdot b_1)\sigma (a_2^{(0)}b_2\ot c)=
\sigma (b_1\ot b_2^{(-1)}\cdot c_1)\sigma (a\ot b_2^{(0)}c_2),
\end{eqnarray}
for all $a, b, c\in B$. Then we can consider the crossed product
$_{\sigma }B=k\# _{\sigma }B$ as in \cite{zhang}, which is an algebra in
$_H^H{\cal YD}$, and whose multiplication is:
\begin{eqnarray}
&&b\cdot b'=\sigma (b_1\ot b_2^{(-1)}\cdot b'_1)b_2^{(0)}b'_2.\label{b7}
\end{eqnarray}
Since $_{\sigma }B$ is an algebra in $_H^H{\cal YD}$, it is in
particular a left $H$-module algebra, so we can consider the
smash product $_{\sigma }B\# H$. \vspace{2mm}

Let now $\gamma :B\rightarrow k$ be a morphism in $_H^H{\cal
YD}$, that is
\begin{eqnarray}
&&\gamma (h\cdot b)=\varepsilon (h)\gamma (b), \label{b8}\\
&&\gamma (b^{(0)})b^{(-1)}=\gamma (b)1_H, \label{b9}
\end{eqnarray}
for all $h\in H$ and $b\in B$. If $\gamma $ is normalized and
convolution invertible in $_H^H{\cal YD}$, with convolution
inverse $\gamma ^{-1}$ in $_H^H{\cal YD}$, the analogue of the
operator $D^1$ is given in $_H^H{\cal YD}$ by:
$$\begin{array}{lll}
D^1(\gamma )(b\ot b')& =\gamma (b_1)\gamma (b_2^{(-1)}\cdot b'_1)
\gamma ^{-1}(b_2^{(0)}b'_2) & \\
&=\gamma (b_1)\gamma (b'_1)\gamma ^{-1}(b_2b'_2) & \hfill (by
\;\; (\ref{b8}))
\end{array}$$
that is, $D^1$ is given by the same formula as for ordinary Hopf
algebras. For a morphism $\gamma :B\rightarrow k$ in $_H^H{\cal
YD}$, the laziness condition is identical to the usual one:
$\gamma (b_1)b_2=b_1\gamma (b_2)$ for all $b\in B$.

\begin{theorem}\label{main}
Let $(H,B)$ be a Hopf admissible pair. \vspace{-5pt}
\begin{enumerate}
\itemsep 0pt
\item[(i)] For a normalized left 2-cocycle $\sigma :B\ot B\rightarrow
k$ in $_H^H{\cal YD}$ define $\overline{\sigma }:(B\times H)\ot
(B\times H) \rightarrow k$,
\begin{eqnarray}
&&\overline{\sigma }(b\times h, b'\times h')=\sigma (b\ot h\cdot b')
\varepsilon (h'). \label{v1}
\end{eqnarray}
Then $\overline{\sigma }$ is a normalized left 2-cocycle on
$B\times H$ and we have $_{\sigma }B\# H=$$\;_{\overline{\sigma
}}(B\times H)$ as algebras. Moreover, $\overline{\sigma }$ is
unique with this property.

\item[(ii)] If $\sigma $ is convolution invertible in $_H^H{\cal
YD}$, then $\overline{\sigma }$ is convolution invertible, with
inverse
\begin{eqnarray}
&&\overline{\sigma }^{-1}(b\times h, b'\times h')=
\sigma ^{-1}(b\ot h\cdot b')
\varepsilon (h'), \label{v2}
\end{eqnarray}
where $\sigma ^{-1}$ is the convolution inverse of $\sigma $ in
$_H^H{\cal YD}$.

\item[(iii)] If $\sigma $ is lazy in $_H^H{\cal YD}$, then
$\overline{\sigma }$ is lazy.

\item[(iv)] If $\sigma,\tau :B\ot B\rightarrow k$ are lazy 2-cocycles
in $_H^H{\cal YD}$, then $\overline{\sigma *\tau
}=\overline{\sigma } *\overline{\tau}$, hence the map $\sigma
\mapsto \overline{\sigma }$ is a group homomorphism from
$Z^2_L(B)$ to $Z^2_L(B\times H)$.

\item[(v)] If $\gamma :B\rightarrow k$ is a normalized and
convolution invertible morphism in $_H^H{\cal YD}$, define
$\overline{\gamma }:B\times H \rightarrow k$ by
\begin{eqnarray}
&&\overline{\gamma }(b\times h)=\gamma (b)\varepsilon (h). \label{v3}
\end{eqnarray}
Then $\overline{\gamma }$ is normalized and convolution
invertible and $\overline{D^1(\gamma )}=D^1(\overline{\gamma })$.
If $\gamma $ is lazy in $_H^H{\cal YD}$, then $\overline{\gamma
}$ is also lazy.

\item[(vi)] If $\sigma $ is a lazy 2-coboundary for $B$ in $_H^H{\cal
YD}$, then $\overline{\sigma }$ is a lazy 2-coboundary for
$B\times H$, so the group homomorphism $Z^2_L(B)\rightarrow
Z^2_L(B\times H)$, $\sigma \mapsto \overline{\sigma }$,
factorizes to a group homomorphism $H^2_L(B)\rightarrow
H^2_L(B\times H)$.
\end{enumerate}
\end{theorem}
\begin{proof}
(i) It is easy to see that $\overline{\sigma }$ is normalized. We
will prove that the multiplications in $_{\sigma }B\# H$ and
$_{\overline{\sigma }} (B\times H)$ coincide, and from the
associativity of $_{\sigma }B\# H$ will follow automatically that
$\overline{\sigma }$ is a left 2-cocycle on $B\times H$. We
compute:
$$\begin{array}{lll}
(b\# h)(b'\# h') & =b\cdot (h_1\cdot b')\# h_2h' & \vspace{3pt} \\
 &=\sigma (b_1\ot b_2^{(-1)}\cdot (h_1\cdot b')_1)b_2^{(0)}(h_1\cdot b')_2
\# h_2h' & \vspace{3pt} \\
 &=\sigma (b_1\ot b_2^{(-1)}h_1\cdot b'_1)b_2^{(0)}(h_2\cdot b'_2)
\# h_3h' & \qquad (by \;\;(\ref{r3})) \vspace{3pt} \\
 & = \sigma (b_1\otimes b_2^{(-1)}h_1\cdot b'_1)
(b_2^{(0)}\times h_2)(b'_2\times h') & \vspace{3pt} \\
 & = \sigma (b_1\otimes b_2^{(-1)}h_1\cdot b'_1)
\varepsilon (b_2^{'(-1)})\varepsilon (h'_1)
(b_2^{(0)}\times h_2)(b_2^{'(0)}\times h'_2) & \vspace{3pt} \\
 & = \overline{\sigma }(b_1\times b_2^{(-1)}h_1, b'_1\times b_2^{'(-1)}h'_1)
(b_2^{(0)}\times h_2)(b_2^{'(0)}\times h'_2) & \vspace{3pt} \\
 & = \overline{\sigma }((b\times h)_1, (b'\times h')_1)(b\times h)_2
(b'\times h')_2 & \vspace{3pt} \\
 & =(b\times h)\cdot _{\overline{\sigma }}(b'\times h').
\end{array}$$
The uniqueness of $\overline{\sigma }$ follows easily by applying
$\varepsilon _B\ot \varepsilon _H$ to the multiplications in
$_{\sigma }B\# H$ and $_{\overline{\sigma }}(B\times
H)$.\vspace{2mm}

(ii) Follows by a direct computation, using the formula
(\ref{b3}) for the convolution in $_H^H{\cal YD}$. \vspace{2mm}

(iii) We have already seen that
\begin{eqnarray*}
&&\overline{\sigma }((b\times h)_1, (b'\times h')_1)(b\times h)_2
(b'\times h')_2=\sigma (b_1\otimes b_2^{(-1)}h_1\cdot b'_1)b_2^{(0)}
(h_2\cdot b'_2)\times h_3h'.
\end{eqnarray*}
Now we compute:
$$\begin{array}{ll}
\overline{\sigma }((b\times h)_2, (b'\times h')_2)(b\times h)_1
(b'\times h')_1 =  & \vspace{3pt} \\
\hspace{1cm} =\overline{\sigma }(b_2^{(0)}\times h_2,
b_2^{'(0)}\times h'_2)(b_1\times b_2^{(-1)}h_1)(b'_1\times b_2^{'(-1)}h'_1) & \vspace{3pt} \\
\hspace{1cm} =\sigma (b_2^{(0)}\ot h_2\cdot b_2^{'(0)})(b_1\times b_2^{(-1)}h_1)(b'_1\times b_2^{'(-1)}h') & \vspace{3pt} \\
\hspace{1cm} =\sigma (b_2^{(0)}\ot h_3\cdot
b_2^{'(0)})b_1((b_2^{(-1)})_1h_1\cdot b'_1)\times (b_2^{(-1)})_2h_2b_2^{'(-1)}h' & \vspace{3pt} \\
\hspace{1cm} =\sigma (b_2^{(0)^{(0)}}\ot h_3\cdot
b_2^{'(0)})b_1(b_2^{(-1)}h_1\cdot b'_1)\times b_2^{(0)^{(-1)}}h_2b_2^{'(-1)}h' & \vspace{3pt} \\
\hspace{1cm} =\sigma (b_2^{(0)^{(0)}}\ot (h_2\cdot
b'_2)^{(0)})b_1(b_2^{(-1)}h_1 \cdot b'_1)\times
b_2^{(0)^{(-1)}}(h_2\cdot b'_2)^{(-1)}h_3h' & \quad (by \;\;(\ref{r6})) \vspace{3pt} \\
\hspace{1cm} =\sigma (b_2^{(0)}\ot h_2\cdot
b'_2)b_1(b_2^{(-1)}h_1\cdot b'_1) \times h_3h' & \quad (by \;\;(\ref{b5})) \vspace{3pt} \\
\hspace{1cm} =\sigma (b_2^{(0)}\ot (h_1\cdot
b')_2)b_1(b_2^{(-1)}\cdot (h_1\cdot b')_1) \times h_2h' & \quad (by \;\;(\ref{r3})) \vspace{3pt} \\
\hspace{1cm} =\sigma (b_1\ot b_2^{(-1)}\cdot (h_1\cdot
b')_1)b_2^{(0)}(h_1\cdot b')_2\times h_2h' & \quad (by \;\;(\ref{b6})) \vspace{3pt} \\
\hspace{1cm} =\sigma (b_1\ot b_2^{(-1)}h_1\cdot
b'_1)b_2^{(0)}(h_2\cdot b'_2)\times h_3h' & \quad (by
\;\;(\ref{r3}))
\end{array}$$
which proves that $\overline{\sigma }$ is indeed lazy.
\vspace{2mm}

(iv) Using the formula (\ref{b3}) for the convolution in
$_H^H{\cal YD}$,
we compute:
$$\begin{array}{lll}
\overline{(\sigma *\tau )}(b\times h, b'\times h') & =
(\sigma *\tau )(b\ot h\cdot b')\varepsilon (h') & \vspace{3pt} \\
&=\sigma (b_1\ot b_2^{(-1)}\cdot (h\cdot b')_1)\tau (b_2^{(0)}\ot
(h\cdot b')_2)\varepsilon (h') & \vspace{3pt} \\
&=\sigma (b_1\ot b_2^{(-1)}h_1\cdot b'_1)\tau (b_2^{(0)}\ot
h_2\cdot b'_2) \varepsilon (h') & (by \;\;(\ref{r3})) \vspace{3pt} \\
&=\sigma (b_1\ot b_2^{(-1)}h_1\cdot b'_1)\varepsilon (b_2^{'(-1)})
\varepsilon (h'_1)\tau (b_2^{(0)}\ot h_2\cdot b_2^{'(0)})\varepsilon (h'_2) & \vspace{3pt} \\
&= \overline{\sigma }(b_1\times b_2^{(-1)}h_1, b'_1\times
b_2^{'(-1)}h'_1)\overline{\tau }(b_2^{(0)}\times h_2, b_2^{'(0)}\times h'_2)& \vspace{3pt} \\
&=(\overline{\sigma }* \overline{\tau })(b\times h, b'\times h'). & \vspace{3pt}  \\
\end{array}$$

(v) Obviously $\overline{\gamma }$ is normalized, and it is easy
to see that its convolution inverse is given by $\overline{\gamma
}^{-1}(b\times h)= \gamma ^{-1}(b)\varepsilon (h)$, where $\gamma
^{-1}$ is the convolution inverse of $\gamma $ in $_H^H{\cal
YD}$. Now we compute:
$$\begin{array}{lll}
\overline{D^1(\gamma )}(b\times h, b'\times h') & = D^1(\gamma )(b\ot h\cdot b')\varepsilon (h') & \vspace{3pt} \\
&=\gamma (b_1)\gamma ((h\cdot b')_1)\gamma ^{-1}(b_2(h\cdot b')_2)\varepsilon (h') & \vspace{3pt} \\
&=\gamma (b_1)\gamma (h_1\cdot b'_1)\gamma ^{-1}(b_2(h_2\cdot
b'_2))\varepsilon (h') & \hspace{-3mm} (by \;\; (\ref{r3})) \vspace{3pt} \\
&=\gamma (b_1)\gamma (b'_1)\gamma ^{-1}(b_2(h\cdot
b'_2))\varepsilon (h') & \hspace{-3mm} (by \;\;(\ref{b8})) \vspace{3pt} \\
&=\gamma (b_1)\gamma (b'_1)\overline{\gamma }^{-1}(b_2(h_1\cdot
b'_2)\times h_2h') & \vspace{3pt} \\
&=\gamma (b_1)\gamma (b'_1)\overline{\gamma }^{-1}((b_2\times h)
(b'_2\times h')) & \vspace{3pt} \\
&=\overline{\gamma }(b_1\times
b_2^{(-1)}h_1)\overline{\gamma}(b'_1\times
b_2^{'(-1)}h'_1)\overline{\gamma}^{-1}((b_2^{(0)}\times h_2)(b_2^{'(0)}\times h'_2)) & \vspace{3pt} \\
&=D^1(\overline{\gamma })(b\times h, b'\times h'). &
\end{array}$$
Hence we have indeed $\overline{D^1(\gamma
)}=D^1(\overline{\gamma })$. Finally, if $\gamma $ is lazy in
$_H^H{\cal YD}$, then we have:
$$\begin{array}{lll}
\overline{\gamma }((b\times h)_1)(b\times h)_2 & =\gamma
(b_1)(b_2\times h) & \vspace{3pt} \\
 & = \gamma (b_2)(b_1\times h) & \\
 & = \gamma (b_2^{(0)})(b_1\times b_2^{(-1)}h)  & \qquad (by \;\;(\ref{b9})) \vspace{3pt} \\
 & = \overline{\gamma }((b\times h)_2)(b\times h)_1, &
\end{array}$$
where the second equality holds because $\gamma$ is lazy, so
$\overline{\gamma }$ is indeed lazy. \vspace{2mm}

(vi) Follows immediately from (v).
\end{proof}

\begin{remark}{\em
Let $\sigma :B\ot B\rightarrow k$ be a normalized morphism in
$_H^H{\cal YD}$, and define $\overline{\sigma }:(B\times H)\ot
(B\times H) \rightarrow k$ by the formula (\ref{v1}). Then one
can easily prove that, conversely, if $\overline{\sigma }$ is a
left 2-cocycle on $B\times H$, then $\sigma $ is a left 2-cocycle
on $B$ in $_H^H{\cal YD}$. Together with (iii) and (iv) of
Theorem \ref{main}, this proves easily that, if $\sigma $ and
$\tau $ are lazy 2-cocycles on $B$ in $_H^H{\cal YD}$, then
$\sigma *\tau $ is a left 2-cocycle on $B$ in $_H^H{\cal YD}$ (as
we mentioned before, this is quite difficult to prove by hand).}
\end{remark}

\begin{example}{\em Let $H_4$ be Sweedler's Hopf algebra. As an
algebra, $H_4=k\langle G,X\vert G^2=1, X^2=0, GX=-XG\rangle$. The
comultiplication is given by $\Delta(G)=G \otimes G$, $\Delta(X)=1
\otimes X+X\otimes G$, and the antipode is $S(G)=G$ and $S(X)=GX$.
This Hopf algebra is a Radford biproduct of the Hopf algebra
$H=k\Ent_2$ and the Hopf algebra $B=k\langle x\vert x^2=0\rangle$
in $_H^H{\cal YD}$. Let $g$ be the generator of the cyclic group
of order two $\Ent_2$. Then $B$ is a left $H$-module algebra with
the action $g \cdot x=-x$ and a left $H$-comodule (co)algebra with
the coaction $\rho(x)=g \otimes x$. The comultiplication and
counit of $B$ are given by $\Delta(x)=1 \otimes x+x \otimes 1$
and $\varepsilon(x)=0$. The Radford biproduct $B \times H$ is
isomorphic to $H_4$ via $1 \times g \mapsto G, x \times g \mapsto
X$. \par \vspace{2mm}

The group of lazy cocycles of $H_4$ is isomorphic to $k$. Any
lazy cocycle $\sigma$ of $H_4$ is of the form
\begin{eqnarray*}
\begin{tabular}{r|rrrr}
$\sigma$ & $1$ & $G$ & $X$ & $GX$ \\
\hline
$1$  & $1$ & $1$ & $0$ & $0$ \\
$G$  & $1$ & $1$ & $0$ & $0$ \\
$X$  & $0$ & $0$ & $\frac{t}{2}$ & $-\frac{t}{2}$ \\
$GX$ & $0$ & $0$ & $\frac{t}{2}$ & $\frac{t}{2}$
\end{tabular}
\end{eqnarray*}
for some $t \in k$, see \cite[Example 2.1]{bc}. The group
$B_L(H_4)$ is trivial, so $H_L^2(H_4)=Z_L^2(H_4) \cong k$. One
may check that any cocycle $\theta$ in $B$ is of the form
$\theta(1,1)=1, \theta(1,x)=\theta(x,1)=0$ and $\theta(x,x)=s$ for
some $s \in k$. Denote this cocycle by $\theta_s$. It is not
difficult to verify that the map $H_L^2(B) \rightarrow
H_L^2(H_4),\ \theta_{-s/2} \mapsto \overline{\theta_{-s/2}}$ is a
group isomorphism. Indeed this isomorphism holds more generally
for Taft's Hopf algebras $H_{n^2}$ and for the Hopf algebras
$E(n)$. It would be interesting to find some sufficient
conditions in a Radford biproduct $B \times H$  for the map
$H_L^2(B) \rightarrow H_L^2(B \times H$) to be an isomorphism.}
\end{example}

\section{Yetter-Drinfeld data obtained from lazy 2-cocycles}\label{sec5}
\setcounter{equation}{0}

Let $A$ be an $H$-bicomodule algebra, with comodule structures
$A\rightarrow A\ot H$, $a\mapsto a_{<0>}\ot a_{<1>}$ and
$A\rightarrow H\ot A$, $a\mapsto a_{[-1]}\ot a_{[0]}$. We can
consider the {\it Yetter-Drinfeld datum} $(H, A, H)$ as in
\cite{cmz} (the second $H$ is regarded as an $H$-bimodule
coalgebra), and the Yetter-Drinfeld category $_A {\cal YD}(H)^H$,
whose objects are $k$-modules $M$ endowed with a left $A$-action
(denoted by $a\ot m\mapsto a\cdot m$) and a right $H$-coaction
(denoted by $m\mapsto m_{(0)}\ot m_{(1)}$) satisfying the
compatibility condition
\begin{eqnarray}
&&a_{<0>}\cdot m_{(0)}\ot a_{<1>}m_{(1)}=(a_{[0]}\cdot
m)_{(0)}\ot (a_{[0]}\cdot m)_{(1)}a_{[-1]}, \label{1}
\end{eqnarray}
for all $a\in A$ and $m\in M$. \par \smallskip

Let now $\sigma $ be a normalized and convolution invertible lazy
2-cocycle on $H$, and consider the $H$-bicomodule algebra
$H(\sigma )$ and the associated category $_{H(\sigma )} {\cal
YD}(H)^H$; for an object $M$ of this category, the compatibility
(\ref{1}) becomes
\begin{eqnarray}
&&h_1\cdot m_{(0)}\ot h_2m_{(1)}=(h_2\cdot m)_{(0)}\ot (h_2\cdot
m)_{(1)}h_1, \label{2}
\end{eqnarray}
for all $h\in H(\sigma )$ and $m\in M$, which is identical to the
compatibility in the usual Yetter-Drinfeld category $_H {\cal
YD}^H$. Just as for $_H {\cal YD}^H$, it is easy to see that
(\ref{2}) is equivalent to
\begin{eqnarray}
&&(h\cdot m)_{(0)}\ot (h\cdot m)_{(1)}=h_2\cdot m_{(0)}\ot
h_3m_{(1)}S^{-1}(h_1). \label{3}
\end{eqnarray}
Our aim will be to prove that, if $M$ is a finite dimensional
object in $_{H(\sigma )} {\cal YD}(H)^H$, then $End (M)$ and $End
(M)^{op}$ are algebras in $_H {\cal YD}^H$.

\begin{lemma}
\begin{enumerate}
\item[(i)] The map $\Delta $, regarded as a map $\Delta :H\rightarrow
H(\sigma )\ot H(\sigma ^{-1})$, is an algebra map. Consequently,
if $M\in $$\;_{H(\sigma )}{\cal M}$ and $N\in $$\;_{H(\sigma
^{-1})}{\cal M}$ then $M\ot N\in $$\;_H{\cal M}$ with action
$h\cdot (m\ot n)=h_1\cdot m\ot h_2\cdot n$.

\item[(ii)] If moreover $M\in
$$\;_{H(\sigma )} {\cal YD}(H)^H$ and $N\in
$$ \;_{H(\sigma ^{-1})} {\cal YD}(H)^H$, then $M\ot N\in $$\;_H {\cal
YD}^H$, with comodule structure $m\ot n\mapsto (m_{(0)}\ot
n_{(0)})\ot n_{(1)}m_{(1)}$.
\end{enumerate}
\end{lemma}
\begin{proof}
A straightforward computation; note that (i) appears also in
\cite{bc}.
\end{proof}

\begin{proposition}
Let $\sigma$ be a normalized and convolution invertible lazy
2-cocycle on $H$. Let $M\in $$\;_{H(\sigma )} {\cal YD}(H)^H$
finite dimensional. Then:
\begin{enumerate}
\item[(i)] $M^*$ becomes an object in $_{H(\sigma ^{-1})} {\cal
YD}(H)^H$, with the following structures (called ''of type 1''):
\begin{eqnarray}
&&(h\cdot m^*)(m)=m^*(S_1(h)\cdot m), \label{12}\\
&&m^*_{(0)}(m)m^*_{(1)}=m^*(m_{(0)})S^{-1}(m_{(1)}),\label{13}
\end{eqnarray}
for all $h\in H$, $m\in M$, $m^*\in M^*$, where $S_1:H\rightarrow
H$ is the map given by (\ref{4});

\item[(ii)] $M^*$ becomes an object in $_{H(\sigma ^{-1})} {\cal
YD}(H)^H$, with the following structures (called ''of type 2''):
\begin{eqnarray}
&&(h\cdot m^*)(m)=m^*(S_2(h)\cdot m),\label{14} \\
&&m^*_{(0)}(m)m^*_{(1)}=m^*(m_{(0)})S(m_{(1)}),\label{15}
\end{eqnarray}
for all $h\in H$, $m\in M$, $m^*\in M^*$, where $S_2:H\rightarrow
H$ is the map given by (\ref{5}).
\end{enumerate}
If $\sigma $ is trivial, i.e. $M\in $$\;_H{\cal YD}^H$, these are
the usual left and right duals of $M$ in $_H{\cal YD}^H$, see
\cite{cvoz1}.
\end{proposition}
\begin{proof}
We prove only (i), while (ii) is similar and left to the reader
(for (i) we will use (\ref{n+9}), for (ii) one has to use
(\ref{n+5})). First, it is known that $M^*$ is a right
$H$-comodule with structure (\ref{13}), and it is a left
$H(\sigma ^{-1})$-module with structure (\ref{12}) because
$S_1:H(\sigma ^{-1})\rightarrow H(\sigma )$ is an algebra
antihomomorphism. Hence, we only have to prove the
Yetter-Drinfeld compatibility condition (\ref{2}) for $M^*$. We
compute, for all $h\in H$, $m\in M$, $m^*\in M^*$:
$$\begin{array}{ll}
(h_1\cdot m^*_{(0)})(m)h_2m^*_{(1)}=m^*_{(0)}(S_1(h_1)\cdot m)h_2m^*_{(1)} &  \hfill (by \ (\ref{12})) \vspace{3pt} \\
\hspace{2cm} =m^*((S_1(h_1)\cdot m)_{(0)})h_2S^{-1}((S_1(h_1)\cdot m)_{(1)}) &  \hfill (by \ (\ref{13})) \vspace{3pt} \\
\hspace{2cm} =m^*(S_1(h_1)_2\cdot m_{(0)})h_2S^{-1}(S_1(h_1)_3m_{(1)}S^{-1}(S_1(h_1)_1)) &  \hfill (by \ (\ref{3})) \vspace{3pt} \\
\hspace{2cm} =m^*(S(h_2)\cdot m_{(0)})h_4S^{-1}(S(h_1)m_{(1)}S^{-1}(S_1(h_3))) &  \hfill (by \ (\ref{8})) \vspace{3pt} \\
\hspace{2cm} =m^*(S(h_2)\cdot m_{(0)})\sigma ^{-1}(S(h_4),h_5)h_6S^{-1}(h_3)S^{-1}(m_{(1)})h_1 &  \hfill (by \ (\ref{4})) \vspace{3pt} \\
\hspace{2cm} =m^*(S(h_2)\cdot m_{(0)})\sigma ^{-1}(S(h_3),h_4)S^{-1}(m_{(1)})h_1 &  \hfill (by \ (\ref{n+9})) \vspace{3pt} \\
\hspace{2cm} =m^*(S_1(h_2)\cdot m_{(0)})S^{-1}(m_{(1)})h_1 &  \hfill (by\ (\ref{4})) \vspace{3pt} \\
\hspace{2cm} =(h_2\cdot m^*)(m_{(0)})S^{-1}(m_{(1)})h_1 &  \hfill (by\ (\ref{12})) \vspace{3pt} \\
\hspace{2cm} =(h_2\cdot m^*)_{(0)}(m)(h_2\cdot m^*)_{(1)}h_1 &
 \hfill (by \ (\ref{13}))
\end{array}$$
so (\ref{2}) holds.
\end{proof}

\begin{remark} {\em
It was proved in \cite{bc} that, if $M\in $$\;_{H(\sigma )}{\cal M}$, then
$M^*\in $$\;_{H(\sigma ^{-1})}{\cal M}$ with action given by $(h\cdot m^*)(m)=
m^*(\phi _{\sigma ^{-1}}(h)\cdot m)$. Since we have proved that
$\phi _{\sigma ^{-1}}=S_1$, it follows that this action coincides with
(\ref{12}). Also, it should be clear that under our hypothesis that $H$ has
bijective antipode, the monoidal $H^2_L(H)$--category constructed in
\cite{bc} has not only left duality, but also right duality, the right dual
of an object $M$ having the $H$-action given by (\ref{14}).}
\end{remark}

We can prove now the following result, generalizing the
well-known fact (see \cite[Proposition 4.1]{cvoz1}) that if $M$
is a finite dimensional Yetter-Drinfeld module then $End (M)$ and
$End (M)^{op}$ are Yetter-Drinfeld module algebras.

\begin{proposition}\label{end}
Let $\sigma$ be a normalized and convolution invertible lazy
2-cocycle on $H$, and $M\in $$\;_{H(\sigma )}{\cal YD}(H)^H$
finite dimensional. Then: \vspace{-5pt}
\begin{enumerate}
\itemsep 0pt
\item[(i)] $End (M)$ becomes an algebra in
$_H{\cal YD}^H$, with $H$-structures:
\begin{eqnarray}
&&(h\cdot f)(m)=h_1\cdot f(S_1(h_2)\cdot m), \label{16}\\
&&f_{(0)}(m)\ot f_{(1)}=f(m_{(0)})_{(0)}\ot
S^{-1}(m_{(1)})f(m_{(0)})_{(1)}, \label{17}
\end{eqnarray}
for all $h\in H$, $m\in M$ and $f\in End (M)$;

\item[(ii)] $End (M)^{op}$ becomes an algebra in $_H{\cal YD}^H$,
with $H$-structures:
\begin{eqnarray}
&&(h\cdot f)(m)=h_2\cdot f(S_2(h_1)\cdot m), \label{18}\\
&&f_{(0)}(m)\ot f_{(1)}=f(m_{(0)})_{(0)}\ot
f(m_{(0)})_{(1)}S(m_{(1)}). \label{19}
\end{eqnarray}
\end{enumerate}
\end{proposition}
\begin{proof}
(i) Since $M\in $$\;_{H(\sigma )} {\cal YD}(H)^H$ and $M^*\in
$$\;_{H(\sigma ^{-1})} {\cal YD}(H)^H$ with structures of type 1,
$M\ot M^*$ becomes an object in $_H{\cal YD}^H$, and by
transferring its structure to $End (M)$ via the canonical
isomorphism we get exactly (\ref{16}) and (\ref{17}), so $End
(M)\in $$\;_H{\cal YD}^H$. It is clear that $End (M)$ is a right
$H^{op}$-comodule algebra (its comodule and algebra structures do
not depend on $\sigma $), so we only have to prove that $End (M)$
is a left $H$-module algebra. For $h\in H$, $f, f'\in End (M)$
and $m\in M$, we have:
$$\begin{array}{lll}
((h_1\cdot f)(h_2\cdot f'))(m) & =(h_1\cdot f)(h_2\cdot f'(S_1(h_3)\cdot m)) & \vspace{3pt} \\
&=h_1\cdot f(S_1(h_2)\cdot (h_3\cdot f'(S_1(h_4)\cdot m))) & \vspace{3pt} \\
&=h_1\cdot f((S_1(h_2)\cdot _{\sigma }h_3)\cdot f'(S_1(h_4)\cdot m)) & \vspace{3pt} \\
&=h_1\cdot f(f'(S_1(h_2)\cdot m)) & \hfill (by \;\;(\ref{6})) \vspace{3pt} \\
&=(h\cdot (ff'))(m). &
\end{array}$$
The relation $h\cdot id_M=\varepsilon (h)id_M$ follows
immediately from (\ref{6}). \smallskip

(ii) The $H$-structures (\ref{18}) and (\ref{19}) come from the
ones of $M^*\ot M$ via the identification $End (M)=M^*\ot M$,
where $M^*$ is regarded now as an object in $\;_{H(\sigma ^{-1})}
{\cal YD}(H)^H$ with structures of type 2. One can prove that
$End (M)^{op}$ is an algebra in $_H{\cal YD}^H$ by a computation
similar to the one in (i), using this time the relation (\ref{7}).
\end{proof}

Let $\sigma$ be as above and $M\in {}_{H(\sigma )}{\cal M}$, not
necessarily finite dimensional. Define two actions of $H$ on $End
(M)$ by the formulae (\ref{16}) and (\ref{18}). Then one can
check by direct computations that these actions give left
$H$-module structures on $End (M)$, and the computations in the
proof of the previous proposition show that actually $End (M)$ is
a left $H$-module algebra with (\ref{16}) and $End (M)^{op}$ is a
left $H$-module algebra with (\ref{18}). In particular, take
$M=H(\sigma )$ and denote $End (H(\sigma ))$ by $A^{\sigma }$.
Then we recover the result in \cite{cc} that $A^{\sigma }$ is a
left $H$-module algebra, with action $(h\cdot f)(h')= h_1\cdot
_{\sigma }f(S_1(h_2)\cdot _{\sigma }h')$, for all $h, h'\in H$ and
$f\in A^{\sigma }$. We will see below that if $H$ is moreover
finite dimensional then $A^{\sigma }$ becomes an algebra in
$_H{\cal YD}^H$. \smallskip

Assume now that $H$ is finite dimensional and $A$ is an
$H$-bicomodule algebra with notation as before. Then, by results
in \cite{cmz} or \cite{bpv}, the category $_A{\cal YD}(H)^H$ is
isomorphic to the category $_{H^*\bowtie A}{\cal M}$ of left
modules over the diagonal crossed product algebra $H^*\bowtie A$.
If $M\in $$\;_A{\cal YD}(H)^H$, then $M$ becomes a left
$H^*\bowtie A$-module with structure
\begin{eqnarray*}
&&(p\bowtie a)\cdot m=p((a\cdot m)_{(1)})(a\cdot m)_{(0)},
\end{eqnarray*}
for all $p\in H^*$, $a\in A$ and $m\in M$. By taking $A=H(\sigma
)$, where $\sigma $ is a normalized and convolution invertible
lazy 2-cocycle on $H$, we obtain that if $M\in
$$\;_{H(\sigma )} {\cal YD}(H)^H$ then $M\in
$$\;_{H^*\bowtie H(\sigma )}{\cal M}$, with $(p\bowtie h)\cdot
m=p((h\cdot m)_{(1)})(h\cdot m)_{(0)}$. \smallskip

On the other hand, we have seen in Proposition \ref{exte} that we
have $H^*\bowtie H(\sigma )=D(H)(\overline{\sigma })$ as
$D(H)$-bicomodule algebras, where $\overline{\sigma }$ is the
extension of $\sigma $ to $D(H)$ given by the formula
(\ref{1.4}). Hence, we get that $M\in$$\;_{D(H)(\overline{\sigma
})}{\cal M}$. By the previous discussion, we obtain that $End
(M)$ and $End (M)^{op}$ are left $D(H)$-module algebras, with
$D(H)$-actions given respectively by:
\begin{eqnarray}
&&((p\ot h)\cdot f)(m)=(p\ot h)_1\cdot f(\overline{S}_1((p\ot
h)_2)\cdot m), \label{t+1} \\
&&((p\ot h)\cdot f)(m)=(p\ot h)_2\cdot f(\overline{S}_2((p\ot
h)_1)\cdot m), \label{t+2}
\end{eqnarray}
for all $p\in H^*$, $h\in H$, $f\in End (M)$, $m\in M$, where
$\overline{S}_1, \overline{S}_2:D(H)\rightarrow D(H)$ are the
maps given by the formulae (\ref{sb1}), (\ref{sb2}). \smallskip

If $M$ is moreover assumed to be finite dimensional, then by
Proposition \ref{end}, $End (M)$ and $End (M)^{op}$ are algebras
in $_H{\cal YD}^H$. Hence they become left $D(H)$-module algebras,
with $D(H)$-actions on $End (M)$ and $End (M)^{op}$ given by
\begin{eqnarray}
(p\ot h)\cdot f=p((h\cdot f)_{(1)})(h\cdot f)_{(0)}, \label {t+3}
\end{eqnarray}
where $h\cdot f$ is the action (\ref{16}), respectively
(\ref{18}). So in this case we have two $D(H)$-module algebra
structures on $End (M)$ and two on $End (M)^{op}$.
\begin{proposition}
The two $D(H)$-module algebra structures as above on $End (M)$
(respectively on $End (M)^{op}$) coincide, and are given
respectively by:
$$\begin{array}{l}
((p\ot h)\cdot f)(m)=p(S^{-1}(m_{(1)})h_3f(S_1(h_4)\cdot
m_{(0)})_{(1)}S^{-1}(h_1))h_2\cdot f(S_1(h_4)\cdot
m_{(0)})_{(0)},\vspace{3pt} \\
((p\ot h)\cdot f)(m)=p(h_4f(S_2(h_1)\cdot
m_{(0)})_{(1)}S^{-1}(h_2)S(m_{(1)}))h_3\cdot f(S_2(h_1)\cdot
m_{(0)})_{(0)},
\end{array}$$
for all $p\in H^*$, $h\in H$, $f\in End (M)$ and $m\in M$.
\end{proposition}
\begin{proof}
We give the proof only for $End (M)$, the one for $End (M)^{op}$
is similar. We compute first the $D(H)$-module structure of $End
(M)$ obtained using $\overline{\sigma }$. We have:
$$\begin{array}{ll}
((p\ot h)\cdot f)(m)=(p_2\ot h_1)\cdot f(\overline{S}_1(p_1\ot h_2)\cdot m) & \vspace{3pt} \\
\hspace{2cm} =(p_2\ot h_1)\cdot f((\varepsilon \ot S_1(h_2))(S^{*-1}(p_1)\ot 1)\cdot m) & \hfill (by \;(\ref{sd1})) \vspace{3pt}\\
\hspace{2cm} =(p_2\ot h_1)\cdot f((\varepsilon \ot S_1(h_2))\cdot ((S^{*-1}(p_1)\ot 1)\cdot m)) & \hfill (by \;(\ref{lulu})) \vspace{3pt}\\
\hspace{2cm} =(p_2\ot h_1)\cdot f((\varepsilon \ot S_1(h_2))\cdot S^{*-1}(p_1)(m_{(1)})m_{(0)}) & \vspace{3pt}\\
\hspace{2cm} =(p_2\ot h_1)\cdot f(p_1(S^{-1}(m_{(1)}))S_1(h_2)\cdot m_{(0)}) & \vspace{3pt}\\
\hspace{2cm} =(p\leftharpoonup S^{-1}(m_{(1)})\ot h_1)\cdot f(S_1(h_2)\cdot m_{(0)}) & \vspace{3pt}\\
\hspace{2cm} = p\leftharpoonup S^{-1}(m_{(1)})((h_1\cdot
f(S_1(h_2)\cdot m_{(0)}))_{(1)})(h_1\cdot
f(S_1(h_2) \cdot m_{(0)}))_{(0)}  & \vspace{3pt} \\
\hspace{2cm} = p(S^{-1}(m_{(1)})h_3f(S_1(h_4)\cdot
m_{(0)})_{(1)}S^{-1}(h_1))h_2\cdot f(S_1(h_4)\cdot m_{(0)})_{(0)}
& \;\;  \hfill (by \;(\ref{3}))
\end{array}$$
We compute now the $D(H)$-module structure of $End (M)$ coming
from $_H{\cal YD}^H$. We have:
$$\begin{array}{lll}
((p\ot h)\cdot f)(m) & =(p((h\cdot f)_{(1)})(h\cdot f)_{(0)})(m) & \vspace{3pt} \\
&=p((h\cdot f)_{(1)})(h\cdot f)_{(0)}(m) & \vspace{3pt} \\
&=p(S^{-1}(m_{(1)})(h\cdot f)(m_{(0)})_{(1)})(h\cdot f)(m_{(0)})_{(0)} & \hspace{-0.5cm} (by \;(\ref{17})) \vspace{3pt} \\
&=p(S^{-1}(m_{(1)})(h_1\cdot f(S_1(h_2)\cdot m_{(0)}))_{(1)})(h_1\cdot f(S_1(h_2)\cdot m_{(0)}))_{(0)} & \vspace{3pt} \\
&=p(S^{-1}(m_{(1)})h_3f(S_1(h_4)\cdot
m_{(0)})_{(1)}S^{-1}(h_1))h_2\cdot f(S_1(h_4)\cdot
m_{(0)})_{(0)}, &
\end{array}$$
so the two structures coincide.
\end{proof}

Let $H$ be of finite dimension and $A$ an
$H$-bicomodule algebra with notation as before. Then one can
check, by direct computation, that $A\in $$\;_A{\cal YD}(H)^H$,
where $A$ is a left $A$-module by the left regular action $a\cdot
b=ab$ for all $a, b\in A$, and $A$ is a right $H$-comodule with
coaction $A\rightarrow A\ot H,\; a\mapsto a_{\{0\}}\ot
a_{\{1\}}S^{-1}(a_{\{-1\}})$ for all $a\in A.$ Hence, if $\sigma$
is a normalized and convolution invertible lazy 2-cocycle on $H$,
by taking $A=H(\sigma )$ we obtain that $H(\sigma )\in
$$\;_{H(\sigma )}{\cal YD}(H)^H$, with $H(\sigma )$-action
$h\cdot l=h\cdot_{\sigma }l$ for all $h,l\in H$, and right
$H$-comodule structure $H(\sigma )\rightarrow H(\sigma )\ot H,\; h
\mapsto h_2\ot h_3S^{-1}(h_1)$ for all $h\in H(\sigma ).$ By
applying all the above to $H(\sigma )\in {}_{H(\sigma )}{\cal
YD}(H)^H$, we obtain that $A^{\sigma }=End(H(\sigma ))$ and
$End(H(\sigma ))^{op}$ are algebras in $_H{\cal YD}^H$.

\begin{proposition}
Let $\sigma$ be a normalized and convolution invertible lazy
2-cocycle on $H$ and $M\in $$\;_{H(\sigma )}{\cal M}$. If $\sigma
$ is a lazy 2-coboundary, then the $H$-module algebra structure
of $End (M)$ given by (\ref{16}) is strongly inner (afforded by
some algebra map $G: H\rightarrow End (M)$). If moreover $H$ is
finite dimensional and $M\in $$\;_{H(\sigma )}{\cal YD}(H)^H$,
then the $D(H)$-module structure of $End (M)$ given by
(\ref{t+1}) is also strongly inner.
\end{proposition}
\begin{proof}
Since $\sigma $ is a lazy 2-coboundary, there exists $\gamma :
H\rightarrow k$ lazy, normalized and convolution invertible such
that $\sigma =D^1(\gamma )$. Then, by \cite{bc}, the map $\varphi
:H(\sigma )\rightarrow H$, $\;\varphi (h)=\gamma
(h_1)h_2$, is an isomorhism of $H$-bicomodule algebras. Define
$F: H(\sigma )\rightarrow End (M)$, $\;F(h)(m)=h\cdot m$,
which is obviously an algebra map. Hence, the map $G:
H\rightarrow End (M)$, $\;G=F\circ \varphi ^{-1}$, is also an
algebra map. Using the laziness of $\gamma $, we can express $F$
as $\;F(h)=G(\varphi (h))=\gamma (h_1)G(h_2)=\gamma (h_2)G(h_1)\;$.
Using (\ref{6}), it is easy to see that $F$ is convolution
invertible with inverse $F^{-1}=F\circ S_1$, so the action
(\ref{16}) is just the inner action afforded by $F$. Hence, we
can write (\ref{16}) as follows:
\begin{eqnarray*}
h\cdot f&=&F(h_1)\circ f\circ F^{-1}(h_2)\\
&=&G(h_1)\gamma (h_2)\circ f\circ \gamma ^{-1}(h_3)G^{-1}(h_4)\\
&=&G(h_1)\circ f\circ G^{-1}(h_2),
\end{eqnarray*}
thus (\ref{16}) is strongly inner, afforded by $G$.\par
\vspace{2mm}

Assume now that $H$ is finite dimensional and $M\in
$$\;_{H(\sigma )}{\cal YD}(H)^H$. Then we know that $M$ becomes a
left $D(H)(\overline{\sigma })$-module, and, due to the embedding
$H^2_L(H)\rightarrow H^2_L(D(H)), \sigma \mapsto \overline{\sigma
}$, since $\sigma $ is a lazy 2-coboundary
for $H$ then $\overline{\sigma }$ is a lazy 2-coboundary for
$D(H)$ (namely, $\overline{\sigma }=D^1(\overline{\gamma })$,
where $\overline{\gamma }: D(H)\rightarrow k$, $\overline{\gamma
}(p\ot h)=p(1)\gamma (h)$). Hence, we can repeat the above proof
for $\overline{\sigma }$ instead of $\sigma $ and $D(H)$ instead
of $H$, and we obtain that the $D(H)$-module structure on $End
(M)$ given by (\ref{t+1}) is also strongly inner.
\end{proof}

We can prove also a partial converse of this result. Recall from
\cite{bc} the normal subgroups $CoInt (H)$ and $CoInn (H)$ of
$Aut_{Hopf} (H)$. If $\gamma \in Reg^1 (H)$, define $ad (\gamma
):H\rightarrow H$ by $ad (\gamma )=\gamma ^{-1}*id_H*\gamma $;
then $ad (\gamma )\in Aut_{Hopf}(H)$ if and only if $D^1(\gamma
)$ is lazy. $CoInt (H)$ is defined as the set of Hopf algebra
automorphisms of $H$ of the type $ad (\gamma )$. It contains the
subgroup
\begin{eqnarray*}
&&CoInn (H)=\{f\in Aut_{Hopf} (H)/ \exists \;\phi \in Alg (H, k)
\;with \; f=(\phi \circ S)*id_H*\phi \}.
\end{eqnarray*}
Suppose that, for a given Hopf algebra $H$, we have $CoInt
(H)=CoInn (H)$, and we have $\sigma \in Z^2_L (H)$ of the form
$\sigma =D^1(\gamma )$, with $\gamma \in Reg^1 (H)$. Then, by
\cite[Lemma 1.12]{bc}, it follows that $\sigma \in B^2_L (H)$,
that is there exists $\chi \in Reg^1_L (H)$ such that $\sigma
=D^1(\chi )$.

\begin{proposition}
Let $\sigma $ be as above and $M\in $$\;_{H(\sigma )}\cal M$
finite dimensional. If the action (\ref{16}) of $H$ on $End (M)$
is strongly inner (afforded by some algebra map $G:H\rightarrow
End (M)$), then there exists $\gamma :H\rightarrow k$ normalized
and convolution invertible such that $\sigma =D^1(\gamma )$. If
moreover we have $CoInn (H)=CoInt (H)$, then $\sigma $ is a lazy
2-coboundary.
\end{proposition}
\begin{proof}
Denote as before $F: H(\sigma )\rightarrow End (M)$,
$\;F(h)(m)=h\cdot m$, which is an algebra map. We have, for all
$h\in H$ and $f\in End (M)$:
$$h\cdot f=F(h_1)\circ f\circ F^{-1}(h_2)=G(h_1)\circ f\circ
G^{-1}(h_2).$$ Hence, if we define $\gamma :H\rightarrow End (M)$
by $\gamma (h)=G^{-1}(h_1)\circ F(h_2)$, we obtain that $\gamma
(h)\circ f=f\circ \gamma (h)$, and since $End (M)$ is a central
algebra and this relation holds for all $f\in End (M)$, it
follows that actually $\gamma $ is a map from $H$ to $k$.
Obviously $\gamma $ is normalized and convolution invertible, so
we only have to prove that $\sigma =D^1({\gamma})$. \par
\vspace{2mm}

First note that, since $G$ is an algebra map, we have
$G^{-1}=G\circ S$, so $G^{-1}$ is an antialgebra map. Also, since
$F: H(\sigma )\rightarrow End (M)$ is an algebra map, we have
$F(hl)=\sigma ^{-1}(h_1, l_1)F(h_2)F(l_2)$ for all $h, l\in H$.
Now we compute:
\begin{eqnarray*}
\gamma (hl)&=&G^{-1}(h_1l_1)F(h_2l_2)\\
&=&G^{-1}(h_1l_1)\sigma ^{-1}(h_2, l_2)F(h_3)F(l_3)\\
&=&\sigma ^{-1}(h_1, l_1)G^{-1}(h_2l_2)F(h_3)F(l_3)
\end{eqnarray*}
because $\sigma ^{-1}$ is lazy. Hence, we have:
$$\begin{array}{ll}
\sigma(h_1, l_1)\gamma (h_2l_2) & =\sigma (h_1,l_1)\sigma^{-1}(h_2,l_2)G^{-1}(l_3)G^{-1}(h_3)F(h_4)F(l_4) \vspace{3pt} \\
&=G^{-1}(l_1)G^{-1}(h_1)F(h_2)F(l_2) \vspace{3pt} \\
&=G^{-1}(l_1)\gamma (h)F(l_2) \vspace{3pt} \\
&=\gamma (h)G^{-1}(l_1)F(l_2) \vspace{3pt} \\
&=\gamma (h)\gamma (l),
\end{array}$$
so we obtain $\sigma(h,l)=\gamma (h_1)\gamma(l_1)
\gamma^{-1}(h_2l_2)=D^1(\gamma )(h, l).$ In general, we do not
know whether $\gamma $ is lazy or whether there exists another
$\chi :H\rightarrow k$ {\it lazy} such that $\sigma =D^1(\chi )$.
However, if $CoInn (H)=CoInt (H)$, then by \cite{bc} such a
$\chi$ exists, so $\sigma $ is a lazy 2-coboundary in this case.
\end{proof}

\section{Lifting projective representations afforded by lazy
2-cocycles}\label{sec6}
\setcounter{equation}{0}

A theorem of Schur asserts that for any finite group $G$ there
exists a finite central extension $C$ such that any projective
representation of $G$ can be lifted to an ordinary representation
of $C$. This theorem has been generalized by Ioana Boca in
\cite{boca}, who proved that any {\it cocommutative} Hopf algebra
$H$ admits a (cocommutative) central extension $B$ such that any
projective representation of $H$ can be lifted to an ordinary
representation of $B$.
The aim of this section is to further generalize her result, by
proving that {\it any} Hopf algebra $H$ admits a central
extension $B$ such that any projective representation of $H$ {\it
afforded by a lazy 2-cocycle} can be lifted to an ordinary
representation of $B$. Our proof follows closely the one of Boca,
so many details will be skipped. The proof will reveal again how
important is the fact that lazy 2-cocycles form a group. \par
\smallskip

If $H$ is a Hopf algebra and $K$ is a Hopf subalgebra of $H$, then
$K^+$ is defined by $K^+=K\cap Ker (\varepsilon )$. If $K$ is a
central Hopf subalgebra of $H$, then $HK^+=K^+H$ and $HK^+$ is a
Hopf ideal of $H$, so $\overline{H}=H/HK^+$ is a Hopf algebra. A
{\it central extension} of $H$ is a Hopf algebra $B$ together
with a central Hopf subalgebra $A$ such that the Hopf algebra
quotient $B/BA^+$ is isomorphic to $H$ (we denote by $\pi $ the
surjection $B\rightarrow H$ with kernel $BA^+$). Recall now from
\cite[Definition 2.2]{boca} the concept of a projective
representation for a Hopf algebra $H$.

\begin{definition}
If $V$ is a vector space, a linear map $T: H\rightarrow End (V)$
is called a projective representation of $H$ if: \vspace{-5pt}
\begin{enumerate}
\itemsep 0pt
\item[(i)] $T$ is convolution invertible;
\item[(ii)] $T(1)=id_V$;
\item[(iii)] $T(h)T(l)=\alpha (h_1, l_1)T(h_2l_2)$ for all $h, l\in H$,
where $\alpha \in Hom (H\ot H, k)$ is convolution invertible.
\end{enumerate}
\end{definition}

It was proved in \cite{boca} that if $T$ is a projective
representation, then $\alpha $ is a normalized (and convolution
invertible) left 2-cocycle and is uniquely determined by $T$ (it
will be called the {\it cocycle of $T$}, or we say that {\it $T$
is afforded by $\alpha $}). Conversely, one can see that, if a map
$T$ as above satisfies $(ii)$ and $(iii)$, where $\alpha $ is a
normalized and convolution invertible left 2-cocycle, then it
also satisfies $(i)$, its convolution inverse being
$T^{-1}(h)=T(S_1(h))$, where $S_1$ is the map defined by
(\ref{4}). Hence, $T$ is a projective representation if
and only if $V$ is a left $_{\alpha}H$-module. Recall now from
\cite[Definition 2.11]{boca} the concept of {\it lifting} of a
projective representation.

\begin{definition}
If $(B, A)$ is a central extension of a Hopf algebra $H$ and
$T:H\rightarrow End (V)$ is a projective representation of $H$,
then we say that $T$ can be lifted to $B$ if there exists an
ordinary representation (algebra map) $X:B\rightarrow End (V)$ and
an element $\gamma \in Reg (B, k)$, with $\gamma (1)=1$, such
that $X=\gamma *(T\circ \pi )$. Such a representation $X$ is
called a lift of $T$.
\end{definition}

\begin{lemma}\label{lema1}
Let $H$ and $A$ be Hopf algebras with $A$ commutative. If
$\;\sigma :H\ot H\rightarrow A$ is a normalized and convolution
invertible left 2-cocycle (with respect to the trivial action of
$H$ on $A$) which is moreover a coalgebra map and is lazy in the
sense that $$\sigma (h_1, l_1)\ot h_2l_2=\sigma (h_2, l_2)\ot
h_1l_1$$ in $A\ot H$, for all $h, l\in H$, then the crossed
product $B=A\# _{\sigma }H$ is a Hopf algebra with: \vspace{-5pt}
\begin{enumerate}
\itemsep 0pt
\item[(1)] $(a\# h)(c\# l)=ac\sigma (h_1, l_1)\# h_2l_2$, for all $a,
c\in A$ and $h, l\in H$;
\item[(2)] $\Delta (a\# h)=(a_1\# h_1)\ot (a_2\# h_2)$;
\item[(3)] $\varepsilon (a\# h)=\varepsilon (a)\varepsilon (h)$;
\item[(4)] $S(a\# h)=(\sigma ^{-1}(S(h_2), h_3)\# S(h_1))(S(a)\# 1)$;
\item[(5)] The map $\pi :B\rightarrow H$, $\pi (a\# h)=\varepsilon
(a)h$ is a Hopf algebra epimorphism with kernel $BA^+$;
\item[(6)] $\;A\simeq A\# 1$ is a central Hopf subalgebra of $B$ and
$$A=B^{co(H)}:=\{b\in B/(id \ot \pi )\Delta (b)=b\ot 1\}.$$
\end{enumerate}
\end{lemma}
\begin{proof}
We only show how to replace the cocommutativity of $H$ in
\cite[Lemma 2.1]{boca}, by the laziness of $\sigma $, the rest of
the proof is identical to the one in \cite{boca}. Namely, one can
compute as in \cite{boca} that
$$\Delta ((a\# h)(c\# l))=(a_1c_1\sigma (h_1,l_1)\# h_3l_3)\ot
(a_2c_2\sigma (h_2, l_2)\# h_4l_4),$$ for all $a, c\in A$ and
$h,l\in H$, using the fact that $\sigma $ is a coalgebra map, and
$$\Delta (a\# h)\Delta (c\# l)=(a_1c_1\sigma (h_1,l_1)\# h_2l_2)\ot
(a_2c_2\sigma (h_3, l_3)\# h_4l_4),$$ and these are equal because,
since $\sigma $ is lazy, we have $\sigma (h_2, l_2)\ot
h_3l_3=\sigma (h_3, l_3)\ot h_2l_2.$
\end{proof}

Denote by $G$ the group $Z^2_L(H)$ of normalized and convolution
invertible lazy 2-cocycles on $H$. Denote by $A$ the finite dual
$(kG)^0$ of the group algebra $kG$, so $A$ is a commutative Hopf
algebra. We can generalize \cite[Lemma 3.1]{boca} as follows.

\begin{lemma}\label{lema2}
Let $H$, $G$, $A$ be as above. Define $\sigma :H\ot H\rightarrow
(kG)^*$ by $\sigma (h, l)(\alpha )=\alpha (h, l)$, for all $h,
l\in H$ and $\alpha \in G$. Then $Im (\sigma )\subseteq A$ and
the corestriction $\sigma :H\ot H\rightarrow A$ is a
coalgebra map and a normalized and convolution invertible lazy
2-cocycle.
\end{lemma}
\begin{proof}
We only prove that $\sigma $ is lazy, the rest of the proof is
identical to the one in \cite{boca}. Namely, we have to prove
that for all $h, l\in H$ we have the equality $\sigma (h_1,
l_1)\ot h_2l_2=\sigma (h_2, l_2)\ot h_1l_1$ in $(kG)^0\ot H$.
This is equivalent to proving that $\sigma (h_1, l_1)(\alpha
)h_2l_2=\sigma (h_2, l_2)(\alpha )h_1l_1$ for all $\alpha \in G$,
that is, $\alpha (h_1, l_1)h_2l_2=\alpha (h_2, l_2)h_1l_1$ for
all $\alpha \in G$, which is obviously true because $G$ consists
exactly of lazy cocycles.
\end{proof}

The following result generalizes \cite[Proposition 2.9]{boca}.

\begin{lemma}
Let $H$ be a Hopf algebra and $T:H\rightarrow End (V)$ a
projective representation afforded by a lazy 2-cocycle $\alpha $
and let $\;u\in Reg (H, k)$ with $u(1)=1$. If $\;W:=u*T$, then $W$
is a projective representation with cocycle $\delta (u)*\alpha $,
where $\delta (u)(h, l)=u(h_1)u(l_1)u^{-1}(h_2l_2)\;$ for all $h,
l\in H$.
\end{lemma}
\begin{proof}
Obviously $W(1)=id_V$; then one computes immediately that
$$\begin{array}{ll}
W(h)W(l) & =u(h_1)u(l_1)T(h_2)T(l_2) \vspace{3pt} \\
&=u(h_1)u(l_1)\alpha (h_2, l_2)T(h_3l_3) \vspace{3pt} \\
&=u(h_1)u(l_1)u^{-1}(h_2l_2)u(h_3l_3)\alpha (h_4, l_4)T(h_5l_5) \vspace{3pt} \\
&=u(h_1)u(l_1)u^{-1}(h_2l_2)\alpha (h_3, l_3)u(h_4l_4)T(h_5l_5) \vspace{3pt} \\
&=(\delta (u)*\alpha)(h_1, l_1)W(h_2l_2)
\end{array}$$
where in the fourth equality we used the fact that $\alpha$ is
lazy.
\end{proof}

The next result generalizes \cite[Proposition 2.12]{boca}.

\begin{proposition}
Let $H$, $A$, $\sigma $, $B$ be as in Lemma \ref{lema1}.
Then:\vspace{-3pt}
\begin{enumerate}
\itemsep 0pt
\item[(i)] If $X$ is an ordinary representation of $B$ such that
$\lambda :=X_{/A}$ is a scalar function, then, if we define
$T(h)=X(1\# h)$, $T$ is a projective representation of $H$
afforded by the lazy 2-cocycle $\lambda \circ \sigma $ and
moreover $X$ is a lift of $T$;

\item[(ii)] If $(T, \alpha )$ is a projective representation of $H$
afforded by the lazy 2-cocycle $\alpha $ and $X$ is a lift of
$T$, then $\lambda :=X_{/A}$ is a scalar function. Moreover, the
lazy 2-cocycles $\alpha $ and $\lambda \circ \sigma $ are
cohomologous (but not necessarily lazy cohomologous);

\item[(iii)] Let $(T, \alpha )$ be a projective representation of $H$
afforded by the lazy 2-cocycle $\alpha $. Then $T$ can be lifted
to $B$ if and only if there exists an algebra map $\lambda
:A\rightarrow k$ such that $\alpha $ is cohomologous to $\lambda
\circ \sigma $ (via a map $u\in Reg (H, k)$ with $u(1)=1$, but
$u$ not necessarily lazy).
\end{enumerate}
\end{proposition}
\begin{proof}
Follows closely the proof in \cite{boca}. The laziness of $\alpha $
is used through the fact that $\delta (u)*\alpha $ can be
written as
$\;(\delta (u)*\alpha )(h, l)=u(h_1)u(l_1)u^{-1}(h_2l_2)\alpha
(h_3,l_3)=u(h_1)u(l_1)\alpha (h_2, l_2)u^{-1}(h_3l_3)$, and
through the fact that one has to use the previous lemma, where
$\alpha $ is supposed to be lazy.
\end{proof}

We can finally obtain the desired result, generalizing
\cite[Theorem 3.2]{boca}.

\begin{theorem}
Let $H$ be a Hopf algebra. Then there exists a central extension
$B$ of $H$ such that any projective representation of $H$
afforded by a lazy 2-cocycle can be lifted to $B$.
\end{theorem}
\begin{proof}
Take as above $G=Z^2_L(H)$, $A=(kG)^0$, $\sigma :H\ot
H\rightarrow (kG)^0$, $\sigma (h, l)(\alpha )=\alpha (h, l)$ for
all $h, l\in H$ and $\alpha \in G$. By Lemma \ref{lema2}, the
hypotheses of Lemma \ref{lema1} are satisfied, so we can consider
the Hopf algebra $B=A\# _{\sigma }H$, a central extension of $H$.
We prove that any projective representation $T$ of $H$ afforded
by a lazy 2-cocycle $\alpha $ can be lifted to $B$. By the
previous proposition, it is enough to find an algebra map $\lambda
:A\rightarrow k$ such that $\alpha $ is cohomologous to $\lambda
\circ \sigma $. As in \cite{boca}, define $\lambda :A\rightarrow
k$ by $\lambda (F)=F(\alpha )$, for all $F\in A=(kG)^0$. Then we
have
$$(\lambda \circ \sigma )(h, l)=\lambda (\sigma (h, l))=\sigma (h,
l)(\alpha )=\alpha (h, l),$$ hence $\alpha =\lambda \circ \sigma
.$ Then, we have, for $F, G\in A=(kG)^0$:
$$\lambda (FG)=(FG)(\alpha )=F(\alpha )G(\alpha )=\lambda
(F)\lambda (G)$$ because $\alpha $ is grouplike in $kG$, and
$\lambda (\varepsilon )=\varepsilon (\alpha )=1$, hence $\lambda
$ is an algebra map.
\end{proof}


\end{document}